\documentclass[12pt,leqno]{article}
\usepackage{amsfonts,latexsym,amscd}    

\newtheorem{theorem}{Theorem}[section]
\newtheorem{definition}[theorem]{Definition}
\newtheorem{proposition}[theorem]{Proposition}
\newtheorem{corollary}[theorem]{Corollary}
\newtheorem{lemma}[theorem]{Lemma}

\newcommand{\beq}{\begin{equation}}
\newcommand{\eeq}{\end{equation}}

\newcommand{\cd}{\ \cdot \ }
\newcommand{\rarrow}{\rightarrow}

\newcommand{\calE}{{\cal{E}}}

\newcommand{\calN}{{\cal{N}}}
\newcommand{\calA}{{\cal{A}}}

\newcommand{\calH}{{\cal{H}}}

\begin{document}

\title{Orthogonal Polynomials of Types $A$ and $B$ \\ and Related
Calogero Models}

\author{Charles F. Dunkl}
\date{}

\maketitle


\begin{abstract}
There are examples of Calogero-Sutherland models associated to the Weyl
groups of type $A$ and $B$. When exchange terms are added to the
Hamiltonians the systems have non-symmetric eigenfunctions, which can be
expressed as products of the ground state with members of a family of
orthogonal polynomials. These polynomials can be defined and studied by
using the differential-difference operators introduced by the author in
Trans. Amer. Math. Soc. 1989 (311), 167-183. After a description of known
results, particularly from the works of Baker and Forrester, and Sahi; there
is a study of polynomials which are invariant or alternating for parabolic
subgroups of the symmetric group. The detailed analysis depends on using two
bases of polynomials, one of which transforms monomially under group actions
and the other one is orthogonal. There are formulas for norms and
point-evaluations which are simplifications of those of Sahi. For any
parabolic subgroup of the symmetric group there is a skew operator on
polynomials which leads to evaluation at $(1,1,\ldots ,1)$ of the quotient
of the unique skew polynomial in a given irreducible subspace by the minimum
alternating polynomial, analogously to a Weyl character formula. The last
section concerns orthogonal polynomials for the type $B$ Weyl group with an
emphasis on the Hermite-type polynomials. These can be expressed by using
the generalized binomial coefficients. A complete basis of eigenfunctions of
Yamamoto's $B_N$ spin Calogero model is obtained by multiplying these
polynomials by the ground state.
\end{abstract}

\section{Introduction}

A Calogero-Sutherland model is an exactly solvable quantum many-body
system in one dimension. There are examples associated to the Weyl
groups of type $A$ and $B$, and by the addition of exchange (reflection)
terms the Hamiltonians have non-symmetric eigenfunctions.  In the two
situations described here, the eigenfunctions are polynomials times the
ground state.

The first example consists of $N$ particles on a circle, with particle $j$
being at angle $\theta_j$, $0 \leq \theta_j < 2 \pi$, parameter $k > 0$;
the Hamiltonian is
\beq
\calH_1 := - \sum^N_{i=1}
\left(\frac{\partial}{\partial\theta_i}\right)^2 + \frac{k}{2} \sum_{1
\leq i < j \leq N} \frac{k-(ij)}{\sin^2(\frac{1}{2}(\theta_i-\theta_j))},
\eeq
where $(ij)$ denotes the transposition (``exchange'') $\theta_i
\leftarrow\!\!\rightarrow \theta_j$.

Under the transformation $x_s = \exp(\theta_s \sqrt{-1})$,
\beq
\calH_1 = \sum^N_{i=1} \left(x_i \, \frac{\partial}{\partial
x_i}\right)^2 + 2k \sum_{1 \leq i < j \leq N} \frac{x_i
x_j}{(x_i-x_j)^2} ((ij)-k).
\eeq
The orthogonal polynomials associated to $\calH_1$ are called the
non-symmetric Jack polynomials; see Baker and Forrester \cite{BF1}, Lapointe
and Vinet \cite{LV2}. In Section 2 this Hamiltonian will be further described.

The second example to be studied is the $B$-type spin Calogero model of
Yamamoto [\cite{Y}, \cite{YT}]; parameters $k, k_1$:
\begin{eqnarray}
\calH_2 & = & - \sum^N_{i=1} \left(\frac{\partial}{\partial
x_i}\right)^2 + \frac{1}{4} \sum^N_{i=1} x_i^2 + \sum^N_{i=1}
\frac{k_1(k_1-\sigma_i)}{x_i^2} \\[4mm]
&& + \ 2k \sum_{1 \leq i < j \leq N}
\left\{\frac{k-\sigma_{ij}}{(x_i-x_j)^2} +
\frac{k-\tau_{ij}}{(x_i+x_j)^2}\right\}, \nonumber
\end{eqnarray}
where $\sigma_i, \sigma_{ij}, \tau_{ij}$ are the reflections in the
hyperoctahedral group $W_N$, defined by $x \sigma_i = (x_1,
\ldots$,$-\stackrel{i}{x_i},\ldots,x_N)$, $x \sigma_{ij} = (\ldots,
\stackrel{i}{x_j}, \ldots, \stackrel{j}{x_i},\ldots)$, $x \tau_{ij} =
(\ldots, -\stackrel{i}{x_j},\ldots,-\stackrel{j}{x_i},\ldots)$. The
coefficient $\frac{1}{4}$ in $\calH_2$ is a coupling constant; it can be
changed by rescaling $x$; this choice is to use the weight function
$\exp(-|x|^2/2)$, as will be seen later. In Section 5 is the description
of a complete orthogonal system of eigenfunctions of $\calH_2$,
consisting of polynomials times the ground state
$$
\prod_{1 \leq i < j \leq N} |x_i^2-x_j^2|^k\prod^N_{i=1} |x_i|^{k_1} \exp(-|x|^2/4).
$$

The technical foundation of this paper is the algebra of
differential-difference (``Dunkl'') operators associated to a reflection
group \cite{D1}. In Section 2, there is an outline of known results for the
type A case, namely the symmetric group $S_N$ acting on $\mathbb{R}^N$
by permutation of coordinates.  This section includes a discussion of
inner products and self-adjoint operators, and orthogonal decompositions.

Section 3 is concerned with polynomials and operators invariant under
parabolic subgroups of $S_N$; these are subgroups which leave intervals
$\{1,2,3,\ldots,\ell_1\}, \{\ell_1+1, \ldots, \ell_1 + \ell_2\}, \ldots$
invariant. Formulas for the norms of invariant polynomials are obtained.

In Section 4, the alternating or skew polynomials and operators are
examined. There is the construction of an important operator associated
to any interval $\{\ell+1, \ell+2,\ldots,\ell+m\}$ which is skew for the
associated symmetric group, and which commutes with the appropriately
transformed version of the Hamiltonian $\calH_1$. Any polynomial which
is skew-symmetric for a parabolic (Young) subgroup of $S_N$ is divisible
by an appropriate minimal alternating polynomial (a product of
discriminants). The skew operator is used to evaluate the ratio at $x =
(1,1,\ldots,1)$, a generalization of the Weyl dimension formula.

Section 5 addresses the type $B$ situation and shows how the type $A$
polynomials can be used to build type B Hermite polynomials, the
eigenfunctions of the transformed $\calH_2$. This results in a complete
set of eigenfunctions with arbitrary parity, that is for any subset $A
\subset \{1,2,\ldots,N\}$ there are eigenfunctions which are odd in
$x_i$, $i \in A$ and even in $x_i$, $i \notin A$. Previously, only the
cases of all even or all odd parity were studied, the so-called
generalized Laguerre polynomials.

\subsubsection*{Notations Used Throughout}

\begin{itemize}
\item
$Z_+ = \{0,1,2,3,\ldots\}, \ \calN_N = Z^N_+$, the set of compositions;
\item
$\calN_N^P$ is the set of partitions with no more than $N$ nonzero
parts; $\calN_N^P = \{\lambda \in \calN_N: \ \lambda_1 \geq \lambda_2
\geq \lambda_3\cdots \geq \lambda_N \geq 0\}$;
\item
for $\alpha \in \calN_N$, $\alpha^+$ denotes the sorting of $\alpha$ to a
partition; the permutation of $\alpha$ lying in $\calN_N^P$;
\item
for $\alpha, \beta \in \calN_N$, the dominance order is defined by
$\alpha \succeq \beta$ if and only if $\sum^j_{i=1} \alpha_i \geq
\sum^j_{i=1} \beta_i$ for $1 \leq j \leq N$; and $\alpha \succ \beta$
means $\alpha \succeq \beta$ and $\alpha \neq \beta$;
\item
for $w \in S_N$, the symmetric group, and $x \in \mathbb{R}^N$, let $xw
\in \mathbb{R}^N$ be defined by $(xw)_i = x_{w(i)}, \ 1 \leq i \leq N$,
and sgn$(w)$ denotes the sign of $w$;
\item
for a function $f$ on $\mathbb{R}^N$, let $(wf)(x) = f(xw), \ x \in \mathbb{R}^N$;
\item
for $\alpha \in \calN_N$, let $x^\alpha := \prod^N_{i=1} x_i^{\alpha_i}$,
then $w(x^\alpha) = x^{w\alpha}$, where $(w\alpha)_i =
\alpha_{w^{-1}(i)}, \ 1 \leq i \leq N$, $ w \in S_N$;
\item
an interval $[\ell +1, \ell+m] := \{j \in \mathbb{Z}: \ \ell+1 \leq j
\leq \ell+m\}$;
\item
for an interval $I$, let $S_I = \{w \in S_N: \ i \notin I \mbox{ implies
} w(i) = i\}$ (that is, $S_I$ is isomorphic to the symmetric group of $I$);
\item
for an interval $I$, let $\sigma_I$ be the longest element in $S_I$,
that is, $\sigma_I(i) = 2 \ell +m + 1-i, \ i \in I = [\ell+1, \ell+m]$;
\item
for an interval $I$, the associated alternating polynomial is $a_I(x) =
\prod \{x_i-x_j: \ i < j \mbox{ and } i,j \in I\}$;
\item
for $\alpha \in \calN_N, \ |\alpha| := \sum^N_{i=1} \alpha_i, \ \alpha!
:= \prod^N_{i=1} (\alpha_i!)$;
\item
for a set $A$, $\# A$ is the cardinality;
\item
for $\lambda \in \calN^P_N$, and $t \in \mathbb{R}$, the hook length
product is $h(\lambda,t) := \prod^N_{i=1} \prod^{\lambda_i}_{j=1}
(\lambda_i -j +t+k(\#\{ s: \ s > i \mbox{ and } j \leq \lambda_s \leq \lambda_i\})$;
\item
for $\alpha \in \calN_N, \ 1 \leq i \leq N, \ \kappa_i(\alpha) =
Nk-k(\#\{s: \ \alpha_s > \alpha_i\}+\#\{s: \ s < i \mbox{ and } \alpha_s
= \alpha_i\})+\alpha_i+1$ (a frequently used eigenvalue associated to $\alpha$);
\item
for an interval $I$, $\alpha \in \calN_N$ satisfies condition ($ \geq,
I$) respectively ($> , I$) if $i,j \in I$ and $i < j$ implies $\alpha_i
\geq \alpha_j$, respectively $\alpha_i > \alpha_j$;
\item
for two linear operators $A, B$ the commutator is $[A,B] := AB-BA$;
\item
for $t \in \mathbb{R}, \ m \in \mathbb{Z}_+$, the shifted factorial is
$(t)_m = \prod^m_{i=1}(t+i-1)$; for $\lambda \in \calN^P_N$ (and
implicit parameter $k$) the generalized shifted factorial is
$(t)_\lambda := \prod^N_{i=1}(t-(i-1)k)_{\lambda_{i}}$;
\item
$1^N = (1,1,\ldots,1) \in \mathbb{R}^N$.
\end{itemize}

\section{Background}
\setcounter{equation}{0}

We review facts about the non-symmetric Jack polynomials expressed in
two different bases, the relation to the operators introduced by
Cherednik, and the inductive calculation of norm formulas by use of
adjacent transpositions.

The symmetric group $S_N$ acts on $\mathbb{R}^N$ by permutation of
coordinates and thus extends to an action on functions
$$
wf(x) := f(xw), \ x \in \mathbb{R}^N, \ w \in S_N.
$$
For a parameter $k \geq 0$, the type A Dunkl operators are defined by
\beq
T_i = \frac{\partial}{\partial x_i} + k \sum_{j \neq i}
\frac{1}{x_i-x_j} (1-(ij)), \ 1 \leq i \leq N.
\eeq
For each partition $\lambda \in \calN^P_N$ (henceforth partitions will
be assumed to have no more than $N$ parts), there is a space $E_\lambda$
of polynomials, invariant and irreducible under the algebra generated by
$\{T_ix_i: \ 1 \leq i \leq N\}$ and $\{w: \ w \in S_N\}$.
This algebra can be considered as a subalgebra of the degenerate double
affine Hecke algebra of type A (the latter acts on Laurent series and
also contains the multiplications by $x_i^{-1}$, $1 \leq i \leq N$,
Cherednik [C], Kakei [K3]).

There is a simple relationship between the rational and trigonometric
differential-difference operators of type A; indeed, let $x_i =
e^{y_i}$, $1 \leq i \leq n$; and suppose $f$ is a linear combination of
$\{e^{y_i}, e^{-y_i}: \ 1 \leq i \leq N\}$, then
$$
(T_ix_i)f = (1+(N-1)k)f + \frac{\partial f}{\partial y_i} + k \sum_{j
\neq \ell} \frac{1}{e^{y_i-y_j}-1} \ (f-(ij)f).
$$
However, there is no corresponding relationship for type $B$ (the weight
function for rational type $B$ on the $N$-torus only involves one parameter).

For $\alpha \in \calN_N$,
\begin{eqnarray*}
T_ix_i x^\alpha & = & (Nk-k\#\{j: \ \alpha_j >
\alpha_i\}+\alpha_i+1)x^\alpha \\[4mm]
&& - \ k \sum_{\alpha_{j}>\alpha_{i}+1}
\sum_{s=0}^{\alpha_j-\alpha_i-2} x^\alpha x_i^{s+1} x_j^{-s-1} \\[4mm]
&& + \ k \sum_{\alpha_{j}\leq\alpha_{i}}
\sum_{s=1}^{\alpha_i-\alpha_j} x^\alpha x_i^{-s} x_j^{s};
\end{eqnarray*}
every $x^\beta$ in the sums satisfies $\beta^+ \prec \alpha^+$ except
the cases $s = \alpha_i-\alpha_j > 0$ which produce $k(ij) x^\alpha$.
When $j > i$ and $\alpha_i > \alpha_j$, $(ij) \alpha \prec \alpha$.

Thus the operator $U_i := T_i x_i-k \sum_{j < i} (ij)$ satisfies the
triangularity property
$$U_i x^\alpha = \kappa_i(\alpha)x^\alpha +
\sum\{A(\beta,\alpha)x^\beta: \, |\beta| = |\alpha|, \, \beta^+ \prec
\alpha^+ \mbox{ or } \alpha^+ = \beta^+ \mbox{ and } \alpha \succ \beta\}.
$$
The type-A Cherednik operator $\xi_i$ (as in [BF3]) is defined by
$$
\xi_i = \left(1 \over k\right) x_i \, \frac{\partial}{\partial x_i}+ x_i
\sum_{j < i} \ \frac{1-(ij)}{x_i-x_j}
+\sum_{j > i} \ \frac{x_j(1-(ij))}{x_i-x_j} +1-i,
$$
and satisfies $\xi_i = \left(1 \over k\right) (U_i-(k(N-1)+1))$. The set
$\{U_i: i = 1,\ldots,N\}$ is commutative (more details below), thus
there is a basis (for polynomials) of simultaneous eigenfunctions,
called non-symmetric Jack polynomials. The notation $E_\alpha(x;1/k)$ is
used for the normalization having one as leading coefficient (of
$x^\alpha$); that is,
$$
E_\alpha(x;1/k) = x^\alpha + \sum\{A'(\beta,\alpha)x^\beta: \, |\beta| =
|\alpha|, \, \beta^+ \prec \alpha^+ \mbox{ or } \beta^+ = \alpha^+
\mbox{ and } \beta \prec \alpha\}
$$
(and the coefficients $A'(\beta,\alpha)$ depend on $k$ and $N$).

In the present paper, we use the dual basis $\{p_\alpha: \, \alpha \in
\calN_N\}$ dfined by the generating function
\beq
F_k(x,y) := \sum_\alpha p_\alpha(x)y^\alpha = \prod^N_{i=1}
\left\{(1-x_iy_i)^{-1} \prod^N_{j=1}
(1-x_iy_j)^{-k}\right\}
\eeq
$(x,y \in \mathbb{R}^N)$.  For $\lambda \in \calN^P_N$,
$\omega_\lambda$ is defined to be the scalar multiple of
$E_\lambda(x;1/k)$ such that
$$
\omega_\lambda = p_\lambda + \sum\{B(\beta,\lambda)p_\beta: \, |\beta| =
|\lambda|, \, \beta^+ \succ \lambda\};
$$
the triangularity property of $B(\beta,\lambda)$ was shown in [D3],
further $B(\beta,\lambda)$ is independent of $N$ in the sense that
$B(\beta,\lambda)$ remains constant when $\beta$ and $\lambda$ are
changed to $(\beta_1,\ldots,\beta_N,0,0,\ldots,0)$,
$(\lambda_1,\ldots,\lambda_N,0,\ldots,0) \in \calN_M$ respectively (and
$M \geq N$). Note that the triangularity for $\{p_\alpha\}$ is in the
opposite direction to that of $\{x^\alpha\}$. Next we define the linear
space $E_\lambda$ as the span of the $S_N$-orbit of $\omega_\lambda$,
with basis $\{\omega_\alpha: \, \alpha^+ = \lambda\}$ where
$\omega_{w\lambda} := w \omega_\lambda$ for $w \in S_N$ (this is well
defined, since $w\lambda = \lambda$ implies $w \omega_\lambda =
\omega_\lambda$ for $w \in S_N$).

The intertwining operator of type $A$ is the unique linear map $V$ on
polynomials which satisfies: $V1 = 1$, $V: \, P_n \rarrow P_n$ for each
$n = 0,1,2,\ldots,$ and $V\left({\partial \over \partial x_i} \,
p\right)(x) = T_i (Vp) (x)$ for $1 \leq i \leq N$, $x \in
\mathbb{R}^N$, each polynomial $p$ (see [D3]). Let $\xi$ be the linear
map on polynomials defined by $\xi: \, p_\alpha \mapsto
x^\alpha/\alpha!$, $\alpha \in \calN_N$ and extended by linearity; then
each $E_\lambda$ is an eigenmanifold for $V\xi$ and $V \xi \omega_\alpha
= ((Nk+1)_{\alpha^+})^{-1} \omega_\alpha$, each $\alpha$.

We will use $\langle f,g\rangle$ to denote inner products (of
polynomials $f,g$) which satisfy two conditions: $T_ix_i$ is
self-adjoint for each $i$, and the inner product is $S_N$-invariant;
that is, $\langle T_ix_i f(x),g(x)\rangle = \langle f(x),T_ix_i
g(x)\rangle$ and $\langle wf,g\rangle = \langle f,w^{-1}g\rangle$, $w
\in S_N$.  The irreducibility properties of $E_\lambda$ imply that such
inner products are uniquely determined up to a constant on each $E_\lambda$.

\begin{definition}
For polynomials $f(x) = \sum_\alpha f_\alpha x^\alpha$, $g(x) =
\sum_\alpha g_\alpha x^\alpha$, define the $A$-inner product
$$
\langle f,g\rangle_A := \sum_{\alpha,\beta} f_\alpha g_\beta T^\alpha x^\beta\bigg|_{x=0},
$$
and the $p$-inner product
$$
\langle f,g \rangle_p := \sum_{\alpha,\beta} f_\alpha g_\beta (H^{-1})_{\alpha\beta},
$$
where the matrix $H$ is defined by
$$
F_k(x,y) = \sum_{\alpha,\beta} H_{\alpha\beta} x^\alpha y^\beta.
$$
Alternatively, $\langle x^\alpha, p_\beta(x)\rangle_p = \delta_{\alpha\beta}$.
\end{definition}

Homogeneous polynomials of different degrees are orthogonal in both
inner products.  The $A$-product was introduced in \cite{D2} and shown
to be positive-definite.

\begin{proposition}
The operators $T_ix_i$ are self-adjoint in the $p$- and $A$-inner products.
\end{proposition}

{\em Proof.} The adjoint of multiplication by $x_i$ in the $A$-product
is clearly $T_i$.  For the $p$-product, self-adjointness is equivalent
to
$$
T_i^{(x)} x_i F_k(x,y) = T_i^{(y)} y_i F_k(x,y)
$$
(the superscripts refer to the variables being acted on); but
$$
T_i^{(x)}x_i F_k(x,y)
= F_k(x,y) \left\{1+\frac{(k+1)x_iy_i}{1-x_iy_i} + k \sum_{j \neq i} \frac{1-x_jy_j}{(1-x_iy_j)(1-x_jy_i)}\right\},
$$
which is symmetric under the interchange of $x$ and $y$. \hfill $\Box$

\begin{proposition}
For polynomials $f,g$ and $w \in S_N$, $\langle f,g\rangle_p = \langle
wf,wg\rangle_p$ and $\langle f,g\rangle_A = \langle wf,wg\rangle_A$. In
particular, the transpositions $(ij)$ are self-adjoint.
\end{proposition}

{\em Proof.} The first part follows from the equation $F_k(xw,yw) =
F_k(x,y)$. The second part depends on the transformation properties of
$T_i$, namely $w^{-1} T_i w = T_{w^{-1}(i)}, \ 1 \leq i \leq N, \ w \in
S_N$. \hfill $\Box$

\begin{corollary}
For partitions $\lambda,\mu$ with $\lambda \neq \mu$, $E_\lambda \perp
E_\mu$ in both $p$- and $A$-products.
\end{corollary}

{\em Proof.} The method of (\cite{D3}, Theorem 4.3) used only the
self-adjointness of each $T_ix_i$. \hfill $\Box$
\bigskip

Sahi \cite{Sa} proved this orthogonality for the $p$-product. In [D3] we used
the modification $T_i \rho_i = T_i x_i +k$, where $\rho_i p_\alpha =
p(\alpha_1,\ldots, \alpha_i+1, \ldots)$ ``raising'' operator; see also \cite{D4}.

\begin{proposition}
For $\lambda \in \calN^P_N$, $f,g \in E_\lambda$, $\langle f,g\rangle_A
= (Nk+1)_\lambda \langle f,g\rangle_p$.
\end{proposition}

{\em Proof.} Since $g \in E_\lambda$, $\sum\{p_\alpha T^\alpha g: \
\alpha \in\calN_N, |\alpha| = |\lambda|\} = (Nk+1)_\lambda g$ (formula
for $(V\xi)^{-1}$). That is, $T^\alpha g$ is $(Nk+1)_\lambda$ times the
coefficient of $p_\alpha$ in the expansion of $g$ in the basis
$\{p_\beta\}$.  Let $f = \sum_\alpha f_\alpha x^\alpha$, $g = \sum_\beta
g_\beta p_\beta$, then
\begin{eqnarray*}
\langle f,g\rangle_A & = & \sum_\alpha f_\alpha T^\alpha g = \sum_\alpha
f_\alpha g_\alpha(Nk+1)_\lambda \\[3mm]
& = & (Nk+1)_\lambda \langle f,g\rangle_p. \hspace*{30pt} \Box
\end{eqnarray*}

\begin{corollary} Let $f \in E_\lambda$, then $f(T)^*1 = (Nk+1)_\lambda
f$ where $f(T)^*$ denotes the $p$-adjoint of the operator.
\end{corollary}

{\em Proof.} For any $\mu \in \calN^P_N$, $g \in E_\mu$,
\begin{eqnarray*}
\langle g,f(T)^*1\rangle_p & = & \langle f(T)g,1\rangle_p = \langle
f,g\rangle_A \\[3mm]
& = & \delta_{\mu\lambda}(Nk+1)_\lambda \langle f,g\rangle_p.
\end{eqnarray*}
Since $g$ and $\mu$ are arbitrary, $f(T)^*1 = (Nk+1)_\lambda f$. \hfill $\Box$
\bigskip

In \cite{D3} we showed that
\begin{eqnarray*}
T_i \rho_i \omega_\alpha
& = & (Nk-k \#\{s: \ \alpha_s > \alpha_i\}+\alpha_i +1) \omega_\alpha \\[3mm]
&& + \ k \sum\{(ij) \omega_\alpha: \ \alpha_j > \alpha_i\}, \mbox{ for }
\alpha \in \calN_N, \ 1 \leq i \leq N.
\end{eqnarray*}
Also the commutator $[T_i \rho_i, T_j \rho_j] = k(T_i \rho_i-T_j
\rho_j)(ij)$, for $i \neq j$ (\cite{D3}, Lemma 2.5(iii)).

This leads to the pairwise commuting operators $U_i := T_i \rho_i -k \sum_{j < i}
(ij)$. (These were used by Lapointe and Vinet [\cite{LV1}, \cite{LV2}],
see also Cherednik \cite{C}.)

\begin{definition}
For $\alpha \in \calN_N$, $1 \leq i \leq N$, let
$$
\kappa_i(\alpha) := Nk-k(\#\{s: \ \alpha_s > \alpha_i\} + \#\{s: \ s < i
\mbox{ and } \alpha_s = \alpha_i \})+\alpha_i +1.
$$
\end{definition}

These will appear as eigenvalues of $U_i$, because
\begin{eqnarray*}
U_i \omega_\alpha & = &
\kappa_i(\alpha)\omega_\alpha +k \Sigma\{(ij) \omega_\alpha: \ i < j
\mbox{ and } \alpha_i < \alpha_j\} \\[3mm]
&& - \ k \Sigma \{(ij)\omega_\alpha: \ j < 1
\mbox{ and } \alpha_j < \alpha_i\}.
\end{eqnarray*}
For each partition $\lambda$, the matrix of $U_i$ with respect to the
basis $\{\omega_\alpha: \ \alpha^+ = \lambda\}$ for $E_\lambda$ is
triangular in the dominance ordering (also for the lexicographic order,
a total one). Note if $\alpha \in \calN_N$, and $\alpha_i < \alpha_j, \
i < j$ for some $i, j$, then $(ij) \alpha \succ \alpha$.

The operators $U_i$ satisfy some commutation properties with adjacent
transposition:
\beq
\begin{array}{rl}
{\rm (i)} & [U_i, (j,j+1)] = 0, \ \mbox{if } i < j \mbox{ or }
j+1 < i; \\[2mm]
{\rm (ii)} & (i,i+1) U_i(i,i+1)
= U_{i+1}+k(i,i+1).
\end{array}
\eeq

\begin{theorem}
(\cite{D4}, Section 3)  For each partition $\lambda$ there exists a
unique basis $\{\zeta_\alpha: \ \alpha^+ = \lambda\}$ for $E_\lambda$ satisfying
\begin{itemize}
\item[(1)] $U_i \zeta_\alpha = \kappa_i(\alpha)\zeta_\alpha, \ 1 \leq i
\leq N$;
\item[(2)] $\zeta_\alpha =
\omega_\alpha+\Sigma\{B(\beta,\alpha)\omega_\beta: \ \beta^+ = \lambda
\mbox{ and } \beta \succ \alpha\}$.
\item[(3)] $\langle \zeta_\alpha,\zeta_\beta\rangle = 0$ if $\alpha \neq \beta$.
\end{itemize}
\end{theorem}

\begin{corollary}
$\zeta_\lambda = \omega_\lambda$, and if $\alpha_i = \alpha_{i+1}$, then
$(i,i+1) \zeta_\alpha = \zeta_\alpha$.
\end{corollary}

{\em Proof.} The partition $\lambda$ is the maximum element in
$\{\alpha: \ \alpha^+ = \lambda\}$.

Suppose $\alpha_i = \alpha_{i+1}$, and expand $(i,i+1)\zeta_\alpha$ in the
basis $\{\zeta_\beta: \ \beta^+ = \lambda = \alpha^+\}$. Because $(i,i+1)
\zeta_\alpha$ is an eigenvector for each $U_j$ with $|i-j| > 1$ with
eigenvalue $\kappa_i(\alpha)$, it must be a scalar multiple of
$\zeta_\alpha$. The fact that $(i,i+1)\omega_\alpha = \omega_\alpha$ and (2.3)(ii)
shows the factor is 1. \hfill $\Box$

\begin{proposition}
Suppose $\alpha \in \calN_N$ and $\alpha_i > \alpha_{i+1}$, then
span$\{\zeta_\alpha,\zeta_{\sigma\alpha}\}$ is invariant under $\sigma =
(i,i+1)$, and the matrix of $\sigma$ in this basis is
$$
\left[\begin{array}{cc}
c & 1-c^2 \\[2mm] 1 & -c \end{array}\right] \mbox{ where }
c = \frac{k}{\kappa_i(\alpha)-\kappa_{i+1}(\alpha)}.
$$
\end{proposition}

{\em Proof.} Let $g = \sigma \zeta_\alpha -c \zeta_\alpha$, then $U_j g =
\kappa_j(\alpha)g$ for $j < i$ or $i+1 < j$; this shows $g \in
\mbox{span}\{\zeta_\alpha, \zeta_{\sigma\alpha}\}$. The coefficient of
$\omega_{\sigma\alpha}$ in $g$ is 1, since $\alpha \succ \sigma\alpha$. 
The commutation relation $\sigma U_i \sigma = U_{i+1}+k \sigma$
shows $U_i g = \kappa_{i+1}(\alpha)g$ and $U_{i+1}g =
\kappa_i(\alpha)g$, but $\kappa_{i+1}(\sigma\alpha) = \kappa_i(\alpha)$
and $\kappa_{i+1}(\alpha) = \kappa_i(\sigma\alpha)$, thus $g =
\zeta_{\sigma\alpha}$. Finally, $\sigma g = \zeta_\alpha -c \sigma
\zeta_\alpha = (1-c^2)\zeta_\alpha -cg$. \hfill $\Box$
\bigskip

These equations were found by Sahi \cite{Sa}, see also Baker and
Forrester \cite{BF3}.

In the basis $\{E_\alpha(x;1/k)$, $E_{\sigma\alpha}(x;1/k)\}$ the matrix
of $\sigma$ is $\left[{c \atop 1-c^2}\quad{1 \atop -c}\right]$. By use of the
known evaluations at $x = 1^N$ and the notation of Definitions 3.10 and 3.17,
$$
E_\alpha(x;1/k) = \frac{h(\alpha^+,1)}{h(\alpha^+,k+1)
\calE_+(\alpha)\calE_-(\alpha)} \ \zeta_\alpha (x).
$$

\begin{corollary}
Suppose $\alpha \in \calN_N$ and $\alpha_i > \alpha_{i+1}$, then
$$
\zeta_{\sigma\alpha}(1^N) = (1-c) \zeta_\alpha(1^N),
$$
and
$$
\|\zeta_{\sigma\alpha}\|^2 = (1-c^2) \|\zeta_\alpha\|^2,
$$
for $c = \frac{k}{\kappa_i(\alpha)-\kappa_{i+1}(\alpha)}$ and $\sigma = (i,i+1)$.
\end{corollary}

{\em Proof.} Since $\zeta_{\sigma\alpha} = \sigma\zeta_\alpha-c \zeta_\alpha$,
we have that
\begin{eqnarray*}
\zeta_{\sigma\alpha}(1^N) & = & \sigma \zeta_\alpha(1^N) -c \zeta_\alpha(1^N) \\[2mm]
& = & (1-c)\zeta_\alpha(1^N).
\end{eqnarray*}
Since $\sigma$ is self-adjoint the matrix of $\sigma$ in the orthonormal
basis $\{\zeta_\alpha/\|\zeta_\alpha\|,
\zeta_{\sigma\alpha}/\|\zeta_{\sigma\alpha}\|\}$ must be symmetric, hence
$\|\zeta_{\sigma\alpha}\|^2 = (1-c^2) \|\zeta_\alpha\|^2.$ \hfill $\Box$

\begin{corollary}
In the same notation, let
\beq
f_0 = \zeta_\alpha + \left(\frac{1}{1+c}\right)\zeta_{\sigma\alpha}, \mbox{
and } f_1 = \zeta_\alpha - \left(\frac{1}{1-c}\right) \zeta_{\sigma\alpha},
\eeq
then $\sigma f_0 = f_0$ and $\sigma f_1 = -f_1$.
\end{corollary}

In the next section we derive expressions for the norms $\|\zeta_\alpha
\|^2_p$, $\|\zeta_\alpha\|^2_A$ as a by-product. If $\alpha \in \calN_N$
and $\alpha_i > \alpha_{i+1}$, then
$\kappa_i(\alpha)-\kappa_{i+1}(\alpha) \geq \alpha_i -\alpha_{i+1} +k$,
thus $0 < c < 1$ (when $k > 0$); in fact,
\begin{eqnarray*}
\kappa_i(\alpha)-\kappa_{i+1}(\alpha)
& = & \alpha_i-\alpha_{i+1}+k(1+\#\{s: \ s > i \mbox{ and } \alpha_s =
\alpha_i\} \\[2mm]
&& + \  \#\{s: \ s < i \mbox{ and } \alpha_s = \alpha_{i+1}\} +
\#\{s: \ \alpha_{i+1} < \alpha_s < \alpha_i\}).
\end{eqnarray*}

The relation of the operator $U_i$ to the Hamiltonian $\calH_1$ in
(1.2) is as follows: let
$$
h(x) = \prod_{1 \leq i < j \leq N} |x_i-x_j|^k \prod^N_{i=1} |x_i|^{k(N-1)/2}
$$
(an $S_N$-invariant positively homogeneous function), then
\begin{eqnarray*}
\lefteqn{\hspace*{-30pt}
h(x)(U_i-1-k(N+1)/2)(f(x)/h(x))} \\[4mm]
& := & A_i f(x) = x_i \, \frac{\partial f(x)}{\partial x_i} -k \sum_{j
\neq i} \, \frac{ x_{{\rm max}(i,j)} }{x_i-x_j} \ (ij) f(x),
\end{eqnarray*}
and $\sum^N_{i=1} A^2_i = \calH_1$.  The eigenvalue of $\calH_1$ on the
space $h(x) E_\lambda$ is
$$
\sum^N_{i=1} (\kappa_i(\lambda)-1-k(N+1)/2)^2
= \sum^N_{i=1} \lambda^2_i + k \sum^N_{i=1} (N-2i+1)\lambda_i +k^2 N(N^2-1)/12
$$
(see Baker and Forrester \cite{BF1}, \cite{BF2}, Lapointe and Vinet
\cite{LV2}); the last term is the energy of the ground state. In the
coordinates $x_j = \exp(\sqrt{-1} \, \theta_j)$,
$$
h = 2^k \prod_{i < j} |\sin((\theta_i-\theta_j)/2)|^k.
$$

\section{Subgroup Invariants}
\setcounter{equation}{0}

A parabolic subgroup of $S_N$ is by definition generated by a subset of
$\{(i,i+1): \ 1 \leq i < N\}$. This section concerns subspaces of
$E_\lambda$ invariant under a parabolic subgroup. We start with the
basic structure, an interval and its group of permutations. A typical
interval is denoted by $I$ or $[\ell_1, \ell_2]$ and is defined to be
$\{n \in \mathbb{Z}: \ \ell_1 \leq n \leq \ell_2\}$. The associated
permutation group, denoted by $S_I$, is defined as $\{w \in S_N: \ w(i)
= i \mbox{ for all } i \notin I\}$. Thus $S_I \cong S_m$ where $m = \#
I$.  The parabolic subgroups of $S_N$ are direct products of such groups
corresponding to a collection of disjoint intervals in $[1,N]$.

The technique developed in this section will be used to derive formulas
for the norms of the non-symmetric Jack polynomials $\zeta_\alpha$, as
well as for polynomials with prescribed symmetric or skew-symmetric
properties for parabolic subgroups.

For any $\beta \in \calN_N$, the set $\{w \beta: \ w \in S_I\}$ has a
unique $\succ$-maximal element, which satisfies the following:

\begin{definition}
For an interval $I$, say that a composition $\alpha$ satisfies property
$(\geq, I)$ or $(>,I)$ if $\alpha_i \geq \alpha_j$, respectively $\alpha_i >
\alpha_j$, whenever $i,j \in I$ and $i < j$.
\end{definition}

We deal with the case of one interval first. Let $I = [\ell+1, \ell+m]$
with $1 \leq \ell+1 < \ell+m \leq N$.  The object is to analyze
span$\{\zeta_{w \alpha}: \ w \in S_I\}$ and span$\{w \zeta_\alpha: \ w \in
S_I\}$ for a fixed $\alpha$ satisfying $(\geq,I)$. The structure of
span$\{\zeta_{w\alpha}\}$ mimics that of $E_\lambda$ (with $m$ variables)
with an analogue of $T_i \rho_i$. Part of the motivation for the
following definition is to have commutativity among the operators
associated with disjoint intervals.

\begin{definition}
For a fixed interval $I = [\ell+1, \ell+m]$, for $i \in I$, let
$$
\tau_i := T_i \rho_i-k \sum_{j \leq \ell} (ij).
$$
\end{definition}

Note that $U_i = \tau_i-k \sum_{\ell < j < i}(ij)$, and $w^{-1}
\tau_{w(i)}w = \tau_i$ for $w \in S_I$.

Now fix a composition $\alpha$ satisfying $(\geq,I)$ and let $X =
\mbox{span}\{\zeta_{w\alpha}: \ w \in S_I\}$, a linear space with dim $X =
\#(S_I \alpha)$, a subspace of $E_{\alpha^+}$. Then $X$ is invariant
under $S_I$ because $S_I$ is generated by $\{(i,i+1): \ \ell+1 \leq i <
\ell+m\}$ and Proposition 2.9 applies.

\begin{definition}
For $w \in S_I$ let $g_{w\alpha} := w \zeta_\alpha$. This is well-defined
because $w_1 \alpha = w_2 \alpha$ implies $(w^{-1}_2 w_1) \zeta_\alpha =
\zeta_\alpha$; since the subgroup $\{w \in S_I: \ w\alpha = \alpha\}$ is
generated by adjacent transpositions (from the condition $(\geq, I)$),
Corollary 2.10.
\end{definition}

Note that it is not generally true that $w \beta = \beta$ for some
composition $\beta$ and $w \in S_N$ implies $w \zeta_\beta = \zeta_\beta$,
for example (1,3) $\zeta_{(1,2,1)} \neq \zeta_{(1,2,1)}$.

\begin{proposition}
$X = \mbox{span}\{g_\beta: \ \beta \in S_I \alpha\}$.
\end{proposition}

{\em Proof.} Consider the set $A = \{\gamma: \ \zeta_\gamma \in
\mbox{span}\{g_\beta\}\}$, by Proposition 2.10, if $\gamma \in A$ and
$\gamma_i \neq \gamma_{i+1}$ for $\ell+1 \leq i < \ell+m$, then
$(i,i+1) \gamma \in A$. Also $\alpha \in A$, hence $A = S_I \alpha$.
\hfill $\Box$

\begin{proposition}
For $\beta \in S_I \alpha$, and $i \in I$,
$$
\tau_i g_\beta = \kappa'_i(\beta)g_\beta + k \sum\{(ij)g_\beta: \ j \in
I \mbox{ and } \beta_j > \beta_i\},
$$
where
$$
\kappa_i'(\beta) = Nk-k(\#\{s: \ \beta_s > \beta_i\} + \#\{s: \ s
\leq \ell \mbox{ and } \beta_s = \beta_i\}+\beta_i+1.
$$
\end{proposition}

{\em Proof.} First for $\beta = \alpha$,
\begin{eqnarray}
\tau_i g_\alpha & = & \left(U_i
+k \sum_{\ell < j < i}(ij)\right)g_\alpha \nonumber \\[4mm]
& = & (\kappa_i(\alpha)+k\#\{s: \ \ell < s < i \mbox{ and } \alpha_s =
\alpha_i\})g_\alpha \nonumber \\[4mm]
&& + \ k \sum\{(ij)g_\alpha: \ \ell < j < i \mbox{ and } \alpha_j > \alpha_i\},
\end{eqnarray}
because $g_\alpha = \zeta_\alpha$.  This is the required formula for this
case; the situation $\{j: \ i < j \leq \ell+m \mbox{ and } \alpha_j >
\alpha_i\}$ does not occur.

For an arbitrary $w \in S_I$, let $s = w^{-1}(i), \ \beta = w\alpha$, then
\begin{eqnarray*}
\tau_i g_\beta & = & \tau_i wg_\alpha = w\tau_s g_\alpha \\[3mm]
& = & (Nk-k(\#\{j: \ \alpha_j > \alpha_s\}
+ \#\{j: \ j \leq \ell, \, \alpha_j = \alpha_s\}+\alpha_s+1)wg_\alpha \\[3mm]
&& + \ k \sum\{w(s,j)g_\alpha: \ \ell < j \leq \ell+m, \ \alpha_j > \alpha_s\},
\end{eqnarray*}
but $w(s,j) = (w(s),w(j))w = (i,w(j))w$, $\beta_t = \alpha_{w^{-1}(t)}$
for any $t$, $\beta_i = \alpha_s$.  Thus
$$
\tau_i g_\beta = \kappa'_i(\beta)g_\beta+k \sum\{(i,w(j))g_\beta: \ j
\in I \mbox{ and } \alpha_j > \beta_i\},
$$
and $\alpha_j = \beta_{w(j)}$. \hfill $\Box$
\bigskip

This showed that the structure of the operators $\{\tau_i: \ i \in I\}$
on $X$ is essentially the same as that of $\{T_i \rho_i: \ 1 \leq i \leq
N\}$ on $E_\lambda$.
\bigskip

\begin{proposition}
Suppose $C$ is a linear operator on $X$ and $[C,\tau_i] = 0$ for each $i
\in I$, then $C = c1$, a multiple of the identity.
\end{proposition}

{\em Proof.}  The same proof as (\cite{D3}, Proposition 3.2) for
$E_\lambda$ works, replacing $\{\omega_\beta: \ \beta^+ = \lambda\}$ by
$\{g_\beta: \ \beta \in S_I \alpha\}$. \hfill $\Box$

\begin{proposition}
The operator $U_I := \prod_{i \in I} U_i$ commutes with each
$w \in S_I$; also $[U_I,\tau_i] = 0$ for $i \in I$ and $U_I g = \prod_{i
\in I} \kappa_i(\alpha)g$ for each $g \in X$.
\end{proposition}

{\em Proof.} To show $[U_I,w] = 0$ for $w \in S_I$ it suffices to prove
this for $w = (i,i+1), \ \ell+1 \leq i < \ell+m$. Also $w U_j =
U_j w$ if $j < i$ or $j > i+1$, thus consider
\begin{eqnarray*}
w U_i U_{i+1} w
& = & (U_{i+1} w+k)U_{i+1} w \\
& = & U_{i+1}(U_i w-k)w+k U_{i+1}w \\
& = & U_{i+1} U_i \\
& = & U_i U_{i+1}
\end{eqnarray*}
(by (2.3)). Since $\tau_{\ell+1} = U_{\ell+1}$ we have $[U_I,
\tau_{\ell+1}] = 0$.

For any $w \in S_I$, $\tau_{\ell+1} = w^{-1} \tau_{w(\ell+1)}w$, and
this shows $[U_I, \tau_j] = 0$ for each $j \in I$.

For any basis element $g_{w \alpha}$ of $X$, $U_I g_{w\alpha} = w U_I
g_\alpha = w \prod_{i \in I} \kappa_i(\alpha)g_\alpha$ since $g_\alpha =
\zeta_\alpha$. \hfill $\Box$
\bigskip

The change of basis matrix for $\{g_\beta\}$ to $\{\zeta_\beta\}$ is
triangular; define $B$ by
\beq
\zeta_\beta = \sum_{\gamma \in S_{I}\alpha} B(\gamma,\beta)g_\gamma, \
\mbox{ then } B(\gamma,\gamma) = 1,
\eeq
and $B(\gamma,\beta) = 0$ unless $\gamma \succeq \beta$. The usual proof
(for $\{\omega_\beta\}$ and $\{\zeta_\beta\}$) applies. There is a nice
relationship between $B$ and the Gram matrix for $\{g_\beta\}$.

\begin{definition}
For $\beta, \gamma \in S_I \alpha$, let
$$
H(\beta,\gamma) := \langle g_\beta, g_\gamma\rangle/\|g_\alpha\|^2
$$
(independent of choice of permissible inner product).
\end{definition}

\begin{proposition}
For $w_1, w_2 \in S_I$, and $\beta, \gamma \in S_I \alpha$,
$$
H(w_1\beta,w_2 \gamma) = H(\beta,w_1^{-1} w_2\gamma);
$$
in particular,
$$
H(w_1\alpha,w_2\alpha) = H(\alpha,w_1^{-1} w_2\alpha) =
B^{-1}(\alpha,w_1^{-1} w_2\alpha).
$$
\end{proposition}

{\em Proof.}
The first identity follows from the $S_N$-invariance of the inner
product. For $\beta \in S_I \alpha$,
$$
H(\alpha,\beta) = \langle \zeta_\alpha, \sum_{\gamma \succeq \beta}
B^{-1}(\gamma,\beta) \zeta_\gamma\rangle/\|\zeta_\alpha\|^2 = B^{-1}(\alpha,\beta)
$$
(by orthogonality, $g_\alpha = \zeta_\alpha$). \hfill $\Box$
\bigskip

Sahi \cite{Sa} found a formula for $\|\zeta_\beta\|^2_p$ in terms of a
hook length product associated to the Ferrers diagram of the composition
$\beta$. Here we give an expression whose complexity is roughly the
number of adjacent transpositions needed to transform $\alpha$ to
$\beta$; the upper and lower hook length products for partitions
(Stanley \cite{St}) will also be used eventually.

\begin{definition}
For $\epsilon = + \mbox{ or } -$ (``sign''), an interval $I$, $\beta \in \calN_N$,
let
$$
\calE_\epsilon(\beta,I) := \prod\left\{1+ \frac{\epsilon
k}{\kappa_j(\beta)-\kappa_i(\beta)} : \ \beta_i < \beta_j, \ i < j,
\mbox{ and } i,j \in I\right\}.
$$
\end{definition}

Observe $\calE_\epsilon(\alpha,I) = 1$ ($\alpha$ satisfies $(\geq, I)$).

\begin{lemma}
If $\beta_{i+1} > \beta_i$ for $\ell+1 \leq i < \ell+m$, then 
$$
\calE_\epsilon((i,i+1)\beta,I)/\calE_\epsilon(\beta,I) = 1
+\frac{\epsilon k}{\kappa_i(\beta)-\kappa_{i+1}(\beta)}, \
\epsilon = \pm.
$$
\end{lemma}

{\em Proof.} For $\{i,j\} \subset I$ let $t(\beta;i,j) = 1 +
\frac{\epsilon k}{\kappa_j(\beta)-\kappa_i(\beta)}$ if $\beta_i
< \beta_j$ and $i < j$, else $t(\beta;i,j) = 1$. Recall
$$
\kappa_j(\beta) = Nk-k(\#\{s: \ \beta_s > \beta_j\} + \#\{s: \ s < j, \,
\beta_s = \beta_j\})+\beta_j+1,
$$
each $j$.

Let $\sigma = (i,i+1)$.  Then $t(\sigma \beta;i,j) = t(\beta;i+1,j)$ and
$t(\sigma\beta;i+1,j) = t(\beta;i,j)$ for $j > i+1$; $t(\sigma
\beta;j,i) = t(\beta;j,i+1)$ and $t(\sigma\beta; j,i+1) = t(\beta;j,i)$
for $j < i$. Also $t(\beta;i,i+1) = 1$, and
$$
t(\sigma\beta;i,i+1) = 1 +
\frac{\epsilon k}{\kappa_{i+1}(\sigma\beta)-\kappa_i(\sigma\beta)} = 1
+ \frac{\epsilon k}{\kappa_i(\beta)-\kappa_{i+1}(\beta)}.
$$
The values $t(\beta;j_1,j_2) = t(\sigma\beta;j_1,j_2)$ for indices
$(j_1,j_2)$ not listed above. Since
$\calE_\epsilon(\beta;I) = \prod_{i,j \in I}t(\beta;i,j)$ this shows
$\calE_\epsilon(\sigma\beta;I)/\calE_\epsilon(\beta;I)$ has the
specified value. \hfill $\Box$

\begin{proposition}
Suppose $\beta \in S_I \alpha$, then
$$
\zeta_\beta(1^N) = \calE_-(\beta;I)\zeta_\alpha(1^N),
$$
and
$$
\|\zeta_\beta\|^2 = \calE_+ (\beta;I) \calE_-(\beta;I)\|\zeta_\alpha\|^2.
$$
\end{proposition}

{\em Proof.} Corollary 2.11 showed that
$$
\frac{\|\zeta_{\sigma\beta}\|^2}{\|\zeta_\beta\|^2}
= \frac{\calE_+(\sigma\beta;I) \calE_-(\sigma\beta;I)}{\calE_+(\beta;I) \calE_-(\beta;I)}
$$
for $\beta_i > \beta_{i+1}, \ \ell+1 \leq i < \ell+m$, and $\sigma :=
(i,i+1)$.  The transpositions $(i,i+1)$ generate $S_I$. Similarly,
$\zeta_\beta(1^N)/\calE_-(\beta;I)$ is constant on $S_I \alpha$. \hfill $\Box$
\bigskip

Observe that the case $I = [1,N]$ provides the values
$\|\zeta_\beta\|^2/\|\zeta_\lambda\|^2$ and
$\zeta_\beta(1^N)/\zeta_\lambda(1^N)$ for $\lambda = \beta^+$.

The minimum (for $\succeq$) element in $S_I \alpha$ occurs in several
important formulas. It is denoted by
$$
\alpha ^{R} := (\alpha_1,\ldots,\alpha_\ell,\alpha_{\ell+m},
\alpha_{\ell+m-1},\ldots,\alpha_{\ell+1},\alpha_{\ell+m+1}, \ldots, \alpha_N),
$$
that is, $\alpha^R = \sigma_I \alpha$ where $\sigma_I$ is the longest
element in $S_I$ (the length of a permutation is the minimum number of
adjacent transpositions needed to produce it, as a product).  Thus
$\sigma_I = (\ell+m,\ell+1)(\ell+m-1,\ell+2) \ldots$ and is an
involution, $\sigma^2_I = 1$.

\begin{lemma}
$$
\calE_\epsilon (\alpha^R;I) =
\prod\left\{1+\frac{\epsilon k}{\kappa_i(\alpha)-\kappa_j(\alpha)}: \
\ell+1 \leq i < j \leq \ell+m \mbox{ and } \alpha_i > \alpha_j\right\}.
$$
\end{lemma}

{\em Proof.} When $\alpha_i = \alpha_{i+1}$ for some $i \in I$, then the
list of eigenvalues
$\kappa_{\ell+1}(\alpha^R),\ldots,\kappa_{\ell+m}(\alpha^R)$ is not the
reverse of the list
$\kappa_{\ell+1}(\alpha),\ldots,\kappa_{\ell+m}(\alpha)$, nevertheless
the relative order of pairwise different values is reversed; that is,
$\alpha^R_i > \alpha^R_j$ if and only if $\alpha_{2\ell+m+1-i} >
\alpha_{2 \ell+m+1-j}$ (for $\ell+1 \leq j < i \leq \ell+m$). This
suffices to establish the formula.   \hfill $\Box$
\bigskip

As an example for the case $\alpha_i = \alpha_{i+1}$, take $\ell = 0, \
m = 3$, and $\alpha = (\alpha_1, \alpha_1, \alpha_3)$ with $\alpha_1 >
\alpha_3$; then $\kappa_1(\alpha^R) = \kappa_3(\alpha), \
\kappa_2(\alpha^R) = \kappa_1(\alpha), \ \kappa_3(\alpha^R) = \kappa_2(\alpha)$.

There is a unique $S_I$-invariant in $X$, and if $\alpha$ satisfies
$(>,I)$ there is a unique $S_I$-alternating polynomial in $X$.
Invariance for $S_I$ means $wf = f$, ``alternating'' means $wf =
\mbox{sgn}(w)f$, for all $w \in S_I$.  To establish such properties it
suffices to show $(i,i+1)f = f$ for invariance, $(i,i+1)f = -f$ for
alternating, for $\ell+1 \leq i < \ell+m$.

\begin{definition}
Let
$$
j_{\alpha;I} := \calE_+(\alpha^R;I) \sum_{\beta \in S_{I}\alpha}
\frac{1}{\calE_+(\beta;I)} \ \zeta_\beta;
$$
and
$$
a_{\alpha;I} := \calE_-(\alpha^R;I) \sum_{w \in S_{I}}
\frac{\mbox{sgn}(w)}{\calE_-(w\alpha;I)} \ \zeta_{w\alpha},
$$
provided $\alpha$ satisfies $(>, I)$.
\end{definition}

\begin{theorem}
The polynomials $j_{\alpha;I}$ and $a_{\alpha;I}$ have the following properties:
\begin{itemize}
\item[(1)] $wj_{\alpha;I} = j_{\alpha;I} \mbox{ for } w \in S_I$;
\item[(2)] $\|j_{\alpha;I}\|^2 = (\#S_I \alpha) \calE_+(\alpha^R;I) \|\zeta_\alpha\|^2$;
\item[(3)] $j_{\alpha;I} (1^N) = (\# S_I \alpha) \zeta_\alpha(1^N)$;
\item[(4)] $j_{\alpha;I} = \sum_{\beta \in S_{I} \alpha} \ g_\beta$;
\item[(5)] $w a_{\alpha;I} = (\mbox{sgn } w) a_{\alpha;I} \mbox{ for } w
\in S_I$;
\item[(6)] $\|a_{\alpha;I}\|^2 = (\# S_I) \calE_-(\alpha^R;I) \|\zeta_\alpha\|^2$;
\item[(7)] $a_{\alpha;I} = \sum_{w \in S_{I}}(\mbox{sgn }w) g_{w \alpha}$.
\end{itemize}
\end{theorem}

{\em Proof.} It is clear that the polynomials defined in (4) and (7)
are the only (up to scalar multiple) invariant and alternating elements
of $X$. To show (1) and (5), consider a typical element $f = \sum_{\beta
\in S_{I} \alpha} f_\beta \zeta_\beta$, and fix $i$ with $\ell+1 \leq i <
\ell+m$, and let $\sigma = (i,i+1)$.  Then
\begin{eqnarray*}
f & = & \sum\{f_\beta \zeta_\beta: \ \beta_i = \beta_{i+1}, \ \beta \in
S_I \alpha\} \\[4mm]
&& + \sum\{f_\beta \zeta_\beta + f_{\sigma \beta} \zeta_{\sigma\beta}: \
\beta_i > \beta_{i+1}, \ \beta \in S_I \alpha\}.
\end{eqnarray*}

By Corollary 2.12, $\sigma f = f$ if and only if
$$
\frac{f_{\sigma\beta}}{f_\beta} =
\frac{\calE_+(\beta;I)}{\calE_+(\sigma\beta;I)}
 = \left(1 + \frac{k}{\kappa_i(\beta)-\kappa_{i+1}(\beta)}\right)^{-1},
$$
for each $\beta \in S_I \alpha$ with $\beta_i > \beta_{i+1}$; also
$\sigma f = -f$ if and only if $\beta_i = \beta_{i+1}$ implies $f_\beta
= 0$ and $\beta_i > \beta_{i+1}$ implies
$$
\frac{f_{\sigma\beta}}{f_\beta} =
- \, \frac{\calE_- (\beta;I)}{\calE_-(\sigma\beta;I)}.
$$
Note if $\alpha$ does not satisfy $(>,I)$ then there is no nonzero
alternating element. By the triangularity of $B$, the coefficient
of $g_{\alpha^R}$ in $j_{\alpha,I}$ or $a_{\alpha;I}$ is the same as the
coefficient of $\zeta_{\alpha^R}$; this establishes (4) and (7).

To compute $\|j_{\alpha,I}\|^2$, observe that $\langle g_\beta, \sum_{\gamma
\in S_{I} \alpha} g_\gamma\rangle = \langle g_\alpha,\sum_\gamma
g_\gamma\rangle$ for each $\beta \in S_I \alpha$ by the $S_N$-invariance
of the inner product. Sum this equation over all $\beta \in S_I \alpha$
to obtain
$$
\langle j_{\alpha;I},j_{\alpha;I}\rangle =
(\# S_I \alpha) \langle \zeta_\alpha, j_{\alpha;I}\rangle =
(\# S_I \alpha) \calE_+(\alpha^R;I) \|\zeta_\alpha\|^2
$$
(using the original definition of $j_{\alpha;I}$).

Formula (3) is a direct consequence of (4). The calculation for
$\|a_{\alpha;I}\|^2$ proceeds similarly:
\begin{eqnarray*}
\|a_{\alpha;I}\|^2
& = & \sum_{w_{1} \in S_{I}}\sum_{w_{2} \in S_{I}}
\mbox{ sgn}(w_1) \mbox{sgn}(w_2) \langle
g_{w_{1}\alpha},g_{w_{2}\alpha}\rangle \\[4mm]
& = & (\# S_I) \sum_{w_{2} \in S_{I}}
\mbox{sgn}(w_2) \langle
g_{\alpha},g_{w_{2}\alpha}\rangle \\[4mm]
& = & (\# S_I) \langle \zeta_\alpha, a_{\alpha;I}\rangle \\[4mm]
& = & (\# S_I) \calE_-(\alpha^R;I) \|\zeta_\alpha\|^2
\end{eqnarray*}
(replacing $w_2$ by $w_1w_2$ in the inner sum.) \hfill $\Box$

\begin{corollary}
$$
\sum_{\beta \in S_{I} \alpha} H(\alpha,\beta) = \calE_+(\alpha^R;I),
$$
and if $\alpha$ satisfies $(>,I)$, then
$$
\sum_{w \in S_{I}} \mbox{ sgn}(w)H(\alpha,w \alpha) = \calE_-(\alpha^R;I).
$$
\end{corollary}

{\em Proof.} By definition of $H$,
$$
\sum_{\beta \in S_{I}\alpha}H(\alpha,\beta) = \langle g_\alpha,\sum_\beta
g_\beta\rangle/\|g_\alpha\|^2 = \langle g_\alpha,j_{\alpha;I}\rangle /\|g_\alpha\|^2,
$$
and
$$
\sum_{w \in S_{I}} \mbox{ sgn}(w)H(\alpha,w\alpha) = \langle g_\alpha,a_{\alpha;I}
\rangle/\|g_\alpha\|^2. \hspace*{20pt} \Box
$$
\bigskip

The triangularity argument for extracting the coefficient of
$\zeta_{\alpha^R}$ was used already by Baker, Dunkl, and Forrester
\cite{BDF} in the same context. Earlier, Baker and Forrester \cite{BF2}
considered some special cases of subgroup invariance and relations to
the Jack polynomials.

An analogue of evaluation at $1^N$ for $a_{\alpha;I}$ will be discussed
later. Hook length products are used in the norm calculations.

\begin{definition}
For a partition $\lambda$ and parameter $t$, let
$$
h(\lambda,t) := \prod^{N-m_0}_{i=1} \ \prod^{\lambda_i}_{j=1}
(\lambda_i-j+t+k \#\{s: \ s > i \mbox{ and } j \leq \lambda_s \leq \lambda_i\}),
$$
where $\lambda_i = 0$ for $i > N-m_0$. The special cases $t = 1$, $k$
satisfy $k^{-|\lambda|}h(\lambda,1) = h^*(\lambda)$,
$k^{-|\lambda|}h(\lambda,k) = h_*(\lambda)$, the upper and lower hook
length products for parameter $1/k$ (Stanley \cite{St}), respectively.
\end{definition}

In the next paragraphs, let $I = [1,N]$, and let $\lambda$ be a
partition and $\lambda^R = (\lambda_N, \lambda_{N-1},\ldots,\lambda_1)$.
 We will use known results for the Jack polynomials to determine
$\|\zeta_\lambda\|^2_p$; these are due to Stanley \cite{St}. Sahi first
found $\|\zeta_\lambda\|^2_p$, and recently Baker and Forrester \cite{BF4}
presented a concise self-contained determination of the structural
constants of the Jack polynomials.

We apply the previous results of this section and we suppress the letter
$I$ in $j_{\lambda;I}$ and $\calE_s(\beta;I)$ since $I = [1,N]$.

>From the orthogonality relations on the $N$-torus (more on this later),
it is known that $J_\lambda(x;1/k)$ is a multiple of $j_\lambda$ (where
$J_\lambda$ is the standard Jack polynomial).  Stanley showed
$J_\lambda(1^N;1/k) = (Nk)_\lambda k^{-|\lambda|}$, also in (\cite{D3},
Proposition 4.3) we showed $\zeta_\lambda(1^N) = \omega_\lambda(1^N) =
(Nk+1)_\lambda/h(\lambda,1)$, so for any $\beta \in \calN_N$, by
Proposition 3.12,
\beq
\zeta_\beta(1_N) = \calE_-(\beta)(Nk+1)_{\beta^+}/h(\beta^+,1).
\eeq
Thus by Theorem 3.15(3),
$$
J_\lambda(x;1/k) = \frac{(Nk)_\lambda h(\lambda,1)}{(Nk+1)_\lambda
k^{|\lambda|}(\# S_N \lambda)} \ j_\lambda(x).
$$
Note $\# S_N \lambda = \#\{\beta: \ \beta^+ = \lambda\}$, the dimension
of $E_\lambda$.  The $_1 F_0$ formula for Jack polynomials (Yan
\cite{Yn}, Beerends and Opdam \cite{BO}) asserts
$$
\prod^N_{i=1} (1-x_i)^{-(Nk+1)} = \sum_{\lambda \in \calN^P_N}
\frac{(Nk+1)_\lambda}{k^{|\lambda|}h_*(\lambda)h^*(\lambda)} \ J_\lambda
\left(x;\frac{1}{k}\right).
$$
We convert this to an expression in $j_\lambda$.

\begin{lemma}
$$
h(\lambda,k) = \frac{(Nk)_\lambda \calE_+(\lambda^R)
h(\lambda,k+1)}{(Nk+1)_\lambda(\# S_N \lambda)},
$$
for $\lambda \in \calN^P_N$.
\end{lemma}

{\em Proof.} We will show
$$
\frac{h(\lambda,k)(Nk+1)_\lambda}{h(\lambda,k+1)(Nk)_\lambda} =
\frac{\calE_+ (\lambda^R)}{\# S_N \lambda}.
$$
Denote the factors of $h(\lambda,t)$ by $h(i,j;t) = \lambda_i
-j+t+k\#\{s: \ j \leq \lambda_s \leq \lambda_i\}$, $1 \leq j \leq
\lambda_i$.  Then $h(i,j,k) = h(i,j+1,k+1)$ whenever $j \neq \lambda_s$
for any $s$.  The ratio $h(\lambda,k)/h(\lambda,k+1)$, after
cancellation, is a product of factors like $h(i,\lambda_i,k)/h(i,1,k+1)$
and $h(i,\lambda_s,k)/h(i,\lambda_s+1,k+1)$, for $\lambda_s < \lambda_i$.

We have $h(i,\lambda_i,k) = k(1+\#\{s: \ s > i, \, \lambda_s =
\lambda_i\})$ and
\begin{eqnarray}
h(i,1,k+1)
& = & \lambda_i+k+k \#\{ s: \ s > i, \ 1 \leq \lambda_s \leq \lambda_i\}
\nonumber \\[3mm]
& = & \lambda_i +k(N-m_0-i+1), \nonumber
\end{eqnarray}
where $m_0$ is the number of zero parts of $\lambda$ (as element of
$\calN_N$). For each $j \in \mathbb{Z}_+$, let $m_j = \#\{i: \ \lambda_i =
j\}$, then
$$
\prod^{N-m_0}_{i=1} h(i,\lambda_i,k) = k^{N-m_0} \prod_{j \geq 1} m_j!
$$
Also,
\begin{eqnarray*}
\frac{(Nk+1)_\lambda}{(Nk)_\lambda}
& = & \prod^{N-m_0}_{i=1} \frac{\lambda_i +(N-i+1)k}{k(N-i+1)} \\[4mm]
& = & \frac{m_0!}{k^{N-m_0}N!} \prod^{N-m_0}_{i=1} (\lambda_i+(N-i+1)k).
\end{eqnarray*}
Thus
\begin{eqnarray*}
\lefteqn{\hspace*{-30pt}
\frac{h(\lambda,k)(Nk+1)_\lambda}{h(\lambda,k+1) (Nk)_\lambda}} \\[4mm]
& = & \frac{m_0! \prod_{j \geq 1} m_j!}{N!} \prod^{N-m_0}_{i=1} \prod
\Bigg\{\frac{\lambda_i-\lambda_j+k(1+\#\{s: \ s > i, \, \lambda_j \leq
\lambda_s \leq \lambda_i\})}{\lambda_i-\lambda_j+k(1+\#\{s: \ s > i, \,
\lambda_j < \lambda_s \leq \lambda_i\})} \\[4mm]
&& \hspace*{200pt} : \mbox{ distinct values of } \lambda_j < \lambda_i\Bigg\}.
\end{eqnarray*}

In the latter product $\lambda_j = 0$ is used; it contributes
$$
\frac{\lambda_i +k(N-i+1)}{\lambda_i+k(N-m_0-i+1)}.
$$
Since $\#S_N \lambda = N!/\prod_j m_j!$ it suffices to identify the
remaining product with
$$
\calE_+ (\lambda^R) =
\prod\left\{\frac{\kappa_i(\lambda)-\kappa_j(\lambda)+k}{\kappa_i(\lambda)-\kappa_j(\lambda)}
: \ i < j \mbox{ and } \lambda_i > \lambda_j\right\}.
$$

Let $\mu_1, \mu_2$ be two distinct values in
$(\lambda_1,\lambda_2,\ldots,\lambda_N)$ with $\mu_1 > \mu_2$; supose
$\lambda_{a_{1}} = \lambda_{a_{1}+1} = \cdots = \lambda_{a_{1}+n_{1}-1}
= \mu_1$ and $\lambda_{a_{2}} = \cdots = \lambda_{a_{2}+n_{2}-1} =
\mu_2$ (and no other appearances of $\mu_1, \mu_2$) and $a_1+n_1-1 <
a_2$. This implies $\kappa_t(\lambda) = (N-t+1)k+\mu_1+1$, $a_1 \leq t <
a_1+n_1$ and $\kappa_u(\lambda) = (N-u+1)k+\mu_2+1$, $a_2 \leq u < a_2+n_2$.

The contribution of pair $(\mu_1,\mu_2)$ to the left-hand product is
\begin{eqnarray}
\lefteqn{\hspace*{-30pt}
\prod^{a_1+n_1-1}_{t=a_1} \
\frac{\mu_1-\mu_2+k(a_2+n_2-t)}{\mu_1-\mu_2+k(a_2-t)}} \nonumber \\[4mm]
& = & \prod^{a_1+n_1-1}_{t=a_1} \ \prod^{a_2+n_2-1}_{u=a_2} \
\frac{\mu_1-\mu_2+k(u+1-t)}{\mu_1-\mu_2+k(u-t)} \nonumber \\[4mm]
& = & \prod_{t,u} \
\frac{\kappa_t(\lambda)-\kappa_u(\lambda)+k}{\kappa_t(\lambda)-
\kappa_u(\lambda)}, \nonumber
\end{eqnarray}
which is exactly the contribution of $(\mu_1,\mu_2)$ to
$\calE_+(\lambda^R)$. \hfill $\Box$

\begin{proposition}
$$
\prod^N_{i=1} (1-x_i)^{-(Nk+1)} = \sum_{\lambda \in \calN^P_N} \
\frac{(Nk+1)_\lambda}{h(\lambda,k+1) \calE_+(\lambda^R)} \ j_\lambda(x)
\quad (|x_i| < 1 \mbox{ each } i).
$$
\end{proposition}

{\em Proof.}  This is exactly the $_1 F_0$ series with $h(\lambda,k)$
replaced using the lemma. \hfill $\Box$

\begin{corollary}
For $\lambda \in \calN^P_N$, $\|\zeta_\lambda\|^2_p = \frac{h(\lambda,k+1)}{h(\lambda,1)}$.
\end{corollary}

{\em Proof.}  By the definition of the $p$-inner product,
$$
F_k (x,y) = \sum_{\beta \in \calN_{N}} \ \frac{1}{\|\zeta_\beta\|^2_p} \
\zeta_\beta(x) \zeta_\beta(y).
$$
Set $y = 1_N$, and use $\|\zeta_\beta\|^2_p = \calE_+(\beta)
\calE_-(\beta)\|\zeta_\lambda\|^2_p$,
$$
\zeta_\beta(1^N) = \calE_-(\beta)\zeta_\lambda(1^N) = \calE_-(\beta)(Nk+1)_\lambda/h(\lambda,1),
$$
for $\lambda = \beta^+$. This shows
\begin{eqnarray*}
\prod^N_{i=1} (1-x_i)^{-(Nk+1)}
& = & \sum_{\lambda \in \calN^P_{N}} \
\frac{(Nk+1)_\lambda}{\|\zeta_\lambda\|^2_p h(\lambda,1)} \ \sum_{\beta
\in S_{N} \lambda} \ \frac{\zeta_\beta(x)}{\calE_+ (\beta)} \nonumber \\[4mm]
& = & \sum_{\lambda \in \calN^P_{N}} \
\frac{(Nk+1)_\lambda}{\|\zeta_\lambda\|^2_p h(\lambda,1) \calE_+
(\lambda^R)} \ j_\lambda(x).
\end{eqnarray*}
Match up the coefficients with the $_1 F_0$-expansion. \hfill $\Box$
\bigskip

The value of $\|\zeta_\lambda\|^2_p$ was first obtained by Sahi who used a
recurrence relation (lowering the degree).

There is one more important $S_N$-invariant inner product for which $T_i x_i$ is
self-adjoint, $1 \leq i \leq N$, obtained by considering polynomials as
(analytic) functions on the torus in $\mathbb{C}^N$ (see \cite{D3},
Proposition 4.2).

\begin{definition}
For $f, g$ polynomials with coefficients in $\mathbb{Q}(k)$, let
\beq
\langle f,g\rangle_{\mathbb{T}}
:= \frac{\Gamma(k+1)^N}{(2\pi)^N \Gamma(Nk+1)} \int_{\mathbb{T}^N}
f(x) g^v(x) \left|\prod_{1 \leq j < \ell \leq
N}(x_j-x_\ell)(x_j^{-1}-x_\ell^{-1})\right|^k dm(x),
\eeq
where $g^v(x) := g(x_1^{-1}, x_2^{-1}, \ldots, x_N^{-1})$, $dm(x) = d
\theta_1 \cdots d\theta_N$ and $x_j = \exp(\sqrt{-1} \, \theta_j)$,
$-\pi < \theta_j \leq \pi$, $1 \leq j \leq N$.
\end{definition}

Beerends and Opdam \cite{BO} evaluated $\langle J_\lambda,
J_\lambda\rangle_{\mathbb{T}}$, from which one can deduce: For $\lambda
\in \calN^P_N$, $g \in E_\lambda$,
$$
\|g\|^2_{\mathbb{T}} = \frac{(Nk+1)_\lambda}{((N-1)k+1)_\lambda} \ \|g\|^2_p.
$$
Baker and Forrester \cite{BF3} computed $\|\zeta_\alpha\|^2_{\mathbb{T}}$
with a different method.

\section{Subgroup alternating polynomials}
\setcounter{equation}{0}

Whenever a polynomial is skew-symmetric for a parabolic subgroup of
$S_N$ it is divisible by an appropriate minimal alternating polynomial
(a product of discriminants). The evaluation of the quotient at $x =
1^N$ is a generalization of the Weyl dimension formula. This evaluation
for polynomials in $E_\lambda$ ($\lambda \in \calN^P_N$) will be carried
out in this section, by constructing
for each interval
$I \subset [1,N]$ a skew operator $\psi_I$ with at least these properties:
\begin{itemize}
\item[(1)] $\psi_I w = \mbox{sgn}(w)w \psi_I$, for $w \in S_I$;
\item[(2)] $[\psi_I,U_j] = 0$ for $j \notin I$;
\item[(3)] if $I_1$ is an interval disjoint from $I$, then $[\psi_I,w] =
0$ for $w \in S_{I_{1}}$;
\item[(4)] if $\alpha \in \calN_N$ and satisfies $(>,I)$, then
$\psi_I$ maps span$\{\zeta_{w\alpha:} \ w \in S_I\}$ to itself.
\end{itemize}

Heuristically, one might suspect $\prod_{i < j} (\tau_i-\tau_j)$ works,
but it is not skew because of non-commutativity, and $\prod_{i < j}
(U_i-U_j)$ has the wrong transformation properties (for $(\# I) \geq 3$).

\begin{definition}
For an interval $I$, let $a_I$ denote the minimal alternating polynomial
for $I$, that is, $a_I(x) := \Pi\{x_i-x_j: \ i < j \mbox{ and } i,j
\in I\}$.
\end{definition}

Let $\calA_I$ denote the associated division symmetrizing operator on
polynomials:
$$
\calA_I f(x) := \frac{1}{(\#I)!} \sum_{w \in S_{I}} \mbox{sgn}(w) f(xw)/a_I(x).
$$

For $\alpha$ satisfying $(>,I)$ we will evaluate the functional
$(\calA_I a_{\alpha;I})(1^N)$, in fact, the more general situation for a
collection of disjoint intervals
$((\calA_{I_{1}}\calA_{I_{2}}\cdots)f)(1^N)$, for suitable $f$.

We will impose one more condition on $\psi_I$, which, surprisingly, is
enough to determine the restriction to $X$ uniquely. We will also
construct an operator on $X$, using the Gram matrix $H$, satisfying the
conditions; this will allow the determination of the matrix entries for
$\psi_I$ in the basis $\{g_{w\alpha}: \ w \in S_I\}$.

The aforementioned condition comes from the idea of a ``reversing''
transformation: the requirement that $\psi_I \zeta_{w\alpha}$ is a
simultaneous eigenvector of $\tau_i -k \sum_{i < j \leq \ell+m}(ij)$,
for $\ell < i \leq \ell+m$.  Note that $U_i = \tau_i-k \sum_{\ell < j <
i}(ij)$, so this reverses the interval $I$.

\begin{definition}
For $i \in I$, let $\theta_i := \tau_i - \frac{k}{2} \sum_{j \neq i,
j\in I} (ij)$.
\end{definition}

We show later that if $\psi_I$ is skew and
$[\psi_I,\theta_i] = 0$ for each $i \in I$, then $\psi_I$ has the
reversing property.

Fix $\alpha$ satisfying $(>,I)$, $I = [\ell+1,\ell+m]$ and $X =
\mbox{span}\{g_{w\alpha}: \ w \in S_I\} = \mbox{span}\{\zeta_{w\alpha}: \
w \in S_I\}$.

For a linear transformation $A$ on $X$ we use the matrix notation $A
g_\beta = \sum_{\gamma \in S_{I}\alpha}A(\gamma,\beta)g_\gamma$.

\begin{lemma}
The linear transformation $A$ on $X$ is skew if and only if
$A(w_1\alpha,w_2\alpha) =$ $\mbox{sgn}(w_1)c (w_1^{-1}w_2)$ for
some function $c$ on $S_I$.
\end{lemma}

{\em Proof.} Given the function $c$ define $A$ as indicated (Recall $w_1
\neq w_2$ implies $w_1 \alpha \neq w_2 \alpha$).  The transformation $A$
is skew if and only if $(\mbox{sgn }w_1)w_1A g_{w_{2}\alpha} = Aw_1
g_{w_{2}\alpha}$ for each $w_1, w_2 \in S_I$, that is
\begin{eqnarray}
\mbox{sgn}(w_1) \sum_{w \in S_{I}} A(w \alpha, w_2 \alpha)w_1 g_{w\alpha}
& = & \mbox{sgn}(w_1) \sum_{w \in S_{I}}A(w_1^{-1}
w\alpha,w_2\alpha)g_{w\alpha} \nonumber \\[4mm]
& = & \sum_{w \in S_{I}}A(w\alpha,w_1w_2\alpha) g_{w\alpha}. \nonumber
\end{eqnarray}
Matching up coefficients of $g_{w\alpha}$ shows that $A$ is skew if and
only if
$$
\mbox{sgn}(w_1)A(w_1^{-1} w \alpha,w_2\alpha) =
A(w\alpha,w_1w_2\alpha)
$$
for all $w, w_1, w_2 \in S_I$, consistent with the relation
$$
A(w\alpha,w_1w_2\alpha) = \mbox{sgn}(w)c(w^{-1}w_1w_2). \hspace*{20pt} \Box
$$

\begin{corollary}
With the same hypotheses,
$$
A a_{\alpha;I} = (\sum_{w \in S_{I}}
\mbox{sgn}(w)c(w))j_{\alpha;I}
$$
and
$$
A j_{\alpha;I} = (\sum_{w\in S_{I}}c(w))a_{\alpha;I}.
$$
\end{corollary}

{\em Proof.} By Definition 3.10,
\begin{eqnarray}
A a_{\alpha;I}
& = & \sum_{w_1} \sum_{w_2} \mbox{sgn}(w_2)A(w_1 \alpha_1w_2\alpha)
g_{w_1\alpha} \nonumber \\[4mm]
& = & \sum_{w_1} g_{w_{1}\alpha} \left(\sum_{w_{2}}
\mbox{sgn}(w_2)c(w_1^{-1}w_2)\mbox{sgn}(w_1)\right) \nonumber \\[4mm]
& = & \sum_{w_{1}} g_{w_{1}\alpha}\left(\sum_{w_{3}} \mbox{sgn}(w_3)c(w_3)\right),
\nonumber
\end{eqnarray}
changing the second summation variable $w_2 = w_1w_3$. A similar
argument shows $A j_{\alpha;I} = \left(\sum_w c(w)\right) a_{\alpha;I}$.
\hfill $\Box$
\bigskip

Two more relations apply when $\alpha$ satisfies $(>,I)$:
\beq
\kappa_i(w\alpha) = \kappa_{w^{-1}(i)}(\alpha), \ \mbox{ for } i \in I,
\ w \in S_I;
\eeq
\beq
\calE_\epsilon(w\alpha;I)\calE_\epsilon(\sigma_I w\alpha;I) =
\calE_\epsilon(\alpha^R;I), \mbox{ for } w \in S_I, \ \epsilon = \pm.
\eeq

For the second equation, note for any given $w\alpha$, and $\ell+1 \leq
i < j \leq \ell+m$, the term $1 + \frac{\epsilon
k}{\kappa_i(\alpha)-\kappa_j(\alpha)}$ appears in
$\calE_\epsilon(w\alpha;I)$ if $w(j) < w(i)$, else in
$\calE_\epsilon(\sigma_I w\alpha;I)$; note $(w\alpha)_{w(i)} = \alpha_i$.

\begin{lemma}
Suppose $\alpha$ satisfies $(>,I)$, then
$$
\theta_i g_\beta = \kappa'_i(\beta) g_\beta + \frac{k}{2} \sum_{j\in I,
j \neq i} \mbox{sgn}(\beta_j-\beta_i)(ij)g_\beta,
$$
for $\beta \in S_I \alpha$.
\end{lemma}

{\em Proof.} By Proposition 3.5,
\begin{eqnarray*}
\theta_i g_\beta & = & \kappa'_i(\beta)g_\beta
+ \frac{k}{2} \sum\{(ij) g_\beta: \ j \in I, \ \beta_j > \beta_i\} \\[4mm]
&& - \ \frac{k}{2} \sum\{(ij) g_\beta: \ j \in I, \ \beta_j < \beta_i\}.
\end{eqnarray*}
The case $\beta_i = \beta_j$ cannot occur. \hfill $\Box$
\bigskip

The first important example of a skew operator commuting with each
$\theta_i$ is defined using the Gram matrix for $H$; thus the domain is
just the space $X$ (for a given $\alpha$ satisfying $(>,I)$).

Let $P$ be the operator on $X$ with the matrix $P(w_1\alpha,w_2\alpha) =
\mbox{sgn}(w_1) \delta_{w_{1},w_{2}}$ (``$P$'' suggests parity). Note
that $P^2 = 1$.

\begin{proposition}
The operator on $X$ with the matrix $PH$ is skew and $[PH,\theta_i] = 0$
for $i \in I$.
\end{proposition}

{\em Proof.}
Let $A_i$ be the matrix for $\theta_i$ in the basis $\{g_{w\alpha}: \ w
\in S_I\}$. By the lemma, $A_i(\gamma,\beta) = \kappa'_i(\beta)$ if
$\gamma = \beta$, and $= \frac{k}{2} \mbox{ sgn}(\beta_j-\beta_i)$ if
$\gamma = (i)\beta$, and $=0$ else.  Thus $A_i (\gamma,\beta) = 0$
unless $\gamma = \beta$ or $\gamma = (ij)\beta$ for some $j \in I$, $j
\neq i$.  This shows $PA_iP = A^T_i$ ($T$ for transpose), since
$$
(PA_iP)(w_1\alpha,w_2\alpha) = \mbox{sgn}(w_1)\mbox{sgn}(w_2) A_i(w_1\alpha,w_2\alpha),
$$
thus
$$
(PA_iP)((ij)\beta,\beta) = -A_i((ij)\beta,\beta), \mbox{ for } \beta \in
S_I \alpha.
$$

Also $\theta_i$ is self-adjoint and $H$ is a scalar multiple of the Gram
matrix for $\{g_{w\alpha}: \ w \in S_I\}$, hence $A_i^TH = HA_i$, that
is $PA_i PH = HA_i$ and $A_iPH = PHA_i$.

The operator $PH$ is skew by the lemma (and $H(w_1\alpha,w_2\alpha) =
H(\alpha,w_1^{-1}w_2\alpha), \ w_1,w_2 \linebreak \in S_I$). \hfill $\Box$

\begin{proposition}
Suppose $A$ is a skew operator on $X$, corresponding to the function $c$
on $S_I$, then $[A,\theta_i] = 0$ for each $i \in I$ if and only if
\begin{eqnarray*}
(\kappa_i'(w\alpha)-\kappa'_i(\alpha))c(w) = \frac{k}{2}
\sum_{j\in I,j\neq
i}c((ij)w)(\mbox{sgn}(i-j)+\mbox{sgn}(w^{-1}(j)-w^{-1}(i))), \\[1mm]
\mbox{for each } w \in S_I, \ i \in I.
\end{eqnarray*}
\end{proposition}

{\em Proof.} For a linear transformation $A$ on $X$ the condition
$[A,\theta_i] =0$ is equivalent to
\begin{eqnarray*}
\lefteqn{\hspace*{-30pt}
\kappa'_i(\beta)A(\gamma,\beta) + \frac{k}{2} \sum_{j \neq i} \mbox{sgn}
(\beta_j-\beta_i)A(\gamma,(ij)\beta)} \\[3mm]
& = & \kappa'_i(\gamma)A(\gamma,\beta) + \frac{k}{2} \sum_{j \neq i} \mbox{sgn}
(\gamma_i-\gamma_j)A((ij)\gamma,\beta),
\end{eqnarray*}
($\gamma, \beta \in S_I \alpha, j \in I$). Let $\gamma = w_1 \alpha, \
\beta = w_2 \alpha$ and replace $A(w_1\alpha,w_2\alpha)$ by
$\mbox{sgn}(w_1) c(w_1^{-1}w_2)$.  The equation becomes
\begin{eqnarray*}
\lefteqn{\hspace*{-30pt}
(\kappa'_i(w_2\alpha)-\kappa'_i(w_1\alpha))\mbox{sgn}(w_1)c(w_1^{-1}w_2)}
 \nonumber \\[4mm]
 & = & \frac{k}{2} \mbox{ sgn}(w_1) \sum_{j_{0}\neq i_{0}}
c(w_1^{-1}(w_1(i_0),w_1(j_0))w_2)(\mbox{sgn}(i_0-j_0) \\[4mm]
&& + \ \mbox{sgn}(w_2^{-1}w_1(j_0)-w_2^{-1}w_1(i_0))),
\end{eqnarray*}
where $i_0 = w_1^{-1}(i), \ j_0 = w_1^{-1}(j)$.  Canceling out
$\mbox{sgn}(w_1)$ gives an equation depending only on $w := w_1^{-1}w_2$
and $i_0$, because $\kappa'_i(w_2\alpha) =
\kappa'_{w_{1}(i_{0})}(w_2\alpha) = \kappa'_{i_{0}}(w\alpha)$ and
$\kappa'_i(w_1\alpha) = \kappa'_{i_{0}}(\alpha)$. This is the
equation in the statement. Note for any $w \in S_I$, $i \in I$,
$\kappa'_i(\alpha) = \kappa'_{w(i)}(w\alpha)$, and
$\mbox{sgn}(\alpha_i-\alpha_j) = -\mbox{sgn}(j-i)$, for $i,j \in I$.
\hfill $\Box$

\begin{corollary}
If $A$ is skew and $[A,\theta_i] = 0$ for each $i \in I$, and $k > 0$,
then the values $c(w)$ are uniquely determined for given $c(1)$, $w \in S_I$.
\end{corollary}

{\em Proof.}
A certain subset of the equations in the proposition is extracted. For
any $w \neq 1$ there is a unique $i \in I$ so that $w(j) = j$ for $j <
i$ and $w^{-1}(i) > i$. Specialize the equations to these values of $w,
i$; $\mbox{sgn}(i-j)+\mbox{sgn}(w^{-1}(j)-w^{-1}(i)) \neq 0$ exactly
when $i < j$ and $w^{-1}(i) > w^{-1}(j)$ (by construction  $j < i$
implies $w^{-1}(j) < w^{-1}(i)$).

Thus
$$
(\kappa'_{w^{-1}(i)}(\alpha)-\kappa'_i(\alpha))c(w) = -k \sum\{c((ij)w):
\ w^{-1}(i) > w^{-1}(j), \ i < j \leq \ell +m\}.
$$
The coefficient of $c(w)$ is nonzero, and if $\beta = (ij) w\alpha$ with
$w^{-1}(i) > w^{-1}(j)$, then $\beta \succ w\alpha$. Thus $c(w)$ is
uniquely determined in terms of the values $\{c(w_1): w_1\alpha \succ
w\alpha, \ w_1 \in S_I\}$ when \linebreak $k \neq 0$. \hfill $\Box$
\bigskip

This corollary shows that an operator on polynomials that is skew for
$S_I$ and commutes with $\theta_i, \ i \in I$ is determined on each $X$
($= \mbox{span}\{w \zeta_\alpha: \ w \in S_I\}$, $\alpha$ satisfies
$(>,I)$) as a scalar multiple of $PH$.  We derive some equations which
will be instrumental in computing the multiple.

\begin{proposition} Suppose $A$ is a skew linear operator on $X$ and
$[A,\theta_i] = 0$ for each $i \in I$, then there is a constant $b$ such
that $A = b PH$ and
\begin{itemize}
\item[(1)] $A j_{\alpha;I} = b \calE_+(\alpha^R;I) a_{\alpha;I}$;
\item[(2)] $A a_{\alpha;I} = b \calE_-(\alpha^R;I) j_{\alpha,I}$;
\item[(3)] $A \zeta_{w\alpha} = b \mbox{ sgn}(w)
\calE_+(w\alpha;I)\sigma_I \zeta_{\sigma_{I}w\alpha}, \mbox{ for } w \in S_I$;
\item[(4)] $A^2 = b^2 \calE_+(\alpha^R;I) \calE_-(\alpha^R;I)1$.
\end{itemize}
\end{proposition}

{\em Proof.} The fact that $A$ is a scalar multiple of $PH$ follows from
Corollary 4.8 and Proposition 4.6.  Equations (1) and (2) follow from
Corollary 3.16 and Lemma 4.4. We prove (3) by exhibiting a simultaneous
eigenvector structure for $\{\sigma_I\zeta_{w\alpha}: \ w \in I\}$ and its
relation to $A$.

For $i \in I$, let $U^R_i = \tau_i-k \sum_{i < j \leq \ell+m}(ij)$, then
$\sigma_I U_{\sigma_{I}(i)}\sigma_I = U^R_i$ (note $\sigma_I(i) =
2m+\ell+1-i$).  Further $U^R_i A = AU_i$, indeed $U_i = \theta_i +
\frac{k}{2} \cdot \sum_{j\in I, j \neq i} \mbox{sgn}(j-i)(ij)$ and
$U^R_i = \theta_i - \frac{k}{2} \sum_{j \in I,j\neq i}
\mbox{sgn}(j-i)(ij)$, and $[A,\theta_i] = 0$, $A(ij) = -(ij)A$ by
hypothesis.  For $w \in S_I$, $\sigma_I A \zeta_{w\alpha}$ is an
eigenvector of $U_{\sigma_{I}(i)}$ with eigenvalue $\kappa_i(w\alpha)$,
each $i \in I$, because $\kappa_i(w\alpha)A\zeta_{w\alpha} = A U_i
\zeta_{w\alpha} = U^R_i A \zeta_{w\alpha} = \sigma_I U_{\sigma_{I}(i)}
(\sigma_I A \zeta_{w\alpha})$. Since $\kappa_i(w\alpha) =
\kappa_{\sigma_{I}(i)}(\sigma_I w\alpha)$, this shows there is a
constant $v(w)$ so that $\sigma_I A \zeta_{w\alpha} =
v(w)\zeta_{\sigma_{I}w\alpha}$, each $w \in S_I$. We use (1) to find
$v(w)$; on the one hand
\begin{eqnarray*}
A a_{\alpha;I}
& = & b \calE_- (\alpha^R;I)j_{\alpha;I} \\[4mm]
& = & b \calE_-(\alpha^R;I) \sigma_I j_{\alpha;I} \\[4mm]
& = & b \calE_-(\alpha^R;I) \calE_+(\alpha^R;I) \sigma_I \sum_{w \in S_{I}}
 \frac{1}{\calE_+(w\alpha;I)} \ \zeta_{w\alpha};
 \end{eqnarray*}
on the other hand
$$
A a_{\alpha;I} = \calE_-(\alpha^R;I) \sum_{w\in S_{I}}
\frac{\mbox{sgn}(w)}{\calE_-(w\alpha;I)} \ v(w) \sigma_I \zeta_{I_{\sigma} w\alpha}.
$$
In the first equation, change the summation variable to $\sigma_Iw$, and
match up coefficients  of $\sigma_I \sigma_{w\alpha}$ in the two
equations; this shows
$$
v(w) = b \mbox{ sgn}(w)
\calE_-(w\alpha;I)\calE_+(\alpha^R;I)/\calE_+(\sigma_Iw\alpha;I)
= b \mbox{ sgn}(w) \calE_-(w\alpha;I)\calE_+(w\alpha;I)
$$
(by (4.2)). Finally,
\begin{eqnarray*}
A^2\zeta_{w\alpha} & = & b^2 \mbox{sgn}(\sigma_I) v(w) \sigma_I A
\zeta_{\sigma_{I}w\alpha} \\[4mm]
& - & b^2 \mbox{sgn}(\sigma_I)v(w)v(\sigma_Iw) \zeta_{w\alpha}
\nonumber \\[4mm]
& = & b^2 \calE_+(\alpha^R;I)\calE_-(\alpha^R;I)\zeta_{w \alpha}, \ w \in
S_I. \hspace*{20pt} \Box
\end{eqnarray*}
\bigskip

We construct a skew operator commuting with each $\theta_i$, $i \in I$,
in the algebra generated by $\{\tau_i: \ i \in I\} \cup S_I$, by
induction on the size of the interval.

\begin{definition}
For the interval $I = [\ell+1,\ell+m]$ and $1 \leq s < m$, let $\psi_1
:= 1$,
\begin{eqnarray*}
\tilde{\psi}_{s+1} & := & U_{\ell+1} U_{\ell+2} \cdots U_{\ell+s}
\psi_s; \\[4mm]
\psi_{s+1} & := & \tilde{\psi}_{s+1} - \sum^{\ell+s}_{i=\ell+1}
(i,\ell+s+1) \tilde{\psi}_{s+1} (i,\ell+s+1).
\end{eqnarray*}
Then $\psi_I := \psi_m$.
\end{definition}

\begin{theorem}
The operator $\psi_I$ satisfies
\begin{itemize}
\item[(1)] $\psi_I w = \mbox{sgn}(w) w \psi_I$ for $w \in S_I$;
\item[(2)] $\psi_I w = w \psi_I$ for $w \in S_{[1,\ell]} \times S_{[\ell+m+1,N]}$;
\item[(3)] $[\psi_I, U_j] = 0$ for $j \notin I$;
\item[(4)] $[\psi_I,\theta_i] = 0$, for $i \in I$;
\item[(5)] for any $\alpha$ satisfying $(\geq,I)$, $\mbox{ span}\{w
\zeta_\alpha: \ w \in S_I\}$ is an invariant subspace of $\psi_I$.
\end{itemize}
\end{theorem}

{\em Proof.} Properties (2) and (3) follow immediately from the
defintion; and property (5) is a consequence of (1).  Let (1$_s$) be the
condition $\psi_s(ij) = -(ij) \psi_s$ for $\ell+1 \leq i < j \leq
\ell+s$, and (4$_s$) be
$$
\left[\psi_s, \tau_i - \frac{k}{2} \sum^{\ell+s}_{j=\ell+1,j\neq i}
(ij)\right] = 0.
$$
The case $s = 1$ is trivial. Inductively, suppose (1$_s$) and (4$_s$)
are true. For convenience, let $t = \ell+s+1$: then $(ij)
\tilde{\psi}_{s+1} = - \tilde{\psi}_{s+1}(ij)$ for $\ell+1 \leq i < j
\leq \ell+s$, because $[U_{\ell+1} \cdots U_{\ell+s},(ij)] = 0$ (see
Proposition 3.7). Now evaluate $(ij)\psi_{s+1}(ij)$ using the
definition; each term in the sum is transformed to the negative of the
corresponding term in $\psi_{s+1}$, with the exception of the terms
labeled $i$ and $j$ which are also interchanged; note $(i,j)(i,t) =
(j,t)(i,j)$. It suffices to show 
$(\ell +s,t)\psi_{s+1}(\ell+s,t) = -\psi_{s+1}$; again the terms in the
sum correspond to the negatives, example: $(\ell+s,t)(i,t)
\tilde{\psi}_{s+1}(i,t)(\ell+s,t) = (i,t)(i,\ell+s)
\tilde{\psi}_{s+1}(i,\ell+s)(i,t) = -(i,t)\tilde{\psi}_{s+1}(i,t)$. The
other part is
$(\ell+s,t)(\tilde{\psi}_{s+1}-(\ell+s,t)\tilde{\psi}_{s+1}(\ell+s,t))(\ell+s,t)$.
This shows $(1_{s+1})$.

To show $(4_{s+1})$, let $B_j := \tau_j - \frac{k}{2}
\sum^{\ell+s}_{i=\ell+1,i \neq j}(ij)$, for $\ell+1 \leq j \leq
\ell+s+1$.  Then for $j \leq \ell+s$, $[\tilde{\psi}_{s+1},B_j] = 0$ by
$(4_s)$ and $[U_{\ell+1} \cdots U_{\ell+s},B_j] = 0$; the latter follows
from $\tau_j = (\ell+1,j) U_{\ell+1}(\ell+1,j)$ and $U_{\ell+1} \cdots
U_{\ell+s}$ commutes with $U_{\ell+1}$ and $w \in S_{[\ell+1,\ell+s]}$.

It will suffice to show $[\psi_{s+1},B_t] = 0$ since
$$
(j,t)B_t(j,t) =
\tau_j - \frac{k}{2} \sum^{\ell+s+1}_{i=\ell+1,i \neq j}(ij).
$$

For $\ell+1 \leq i \leq \ell+s$,
\begin{eqnarray*}
[(i,t) \tilde{\psi}_{s+1}(i,t),B_t]
& = & (i,t)[\tilde{\psi}_{s+1},(i,t)B_t(i,t)](i,t) \\[4mm]
& = & (i,t)\left[\tilde{\psi}_{s+1}, B_i - \frac{k}{2} (i,t)\right](i,t)
\\[4mm]
& = & - \, \frac{k}{2} (i,t) [\tilde{\psi}_{s+1},(i,t)](i,t) \\[4mm]
& = & \frac{k}{2} [\tilde{\psi}_{s+1},(it)].
\end{eqnarray*}
Also $[\tilde{\psi}_{s+1},U_t] = 0$ by (3), that is,
$$
\left[\tilde{\psi}_{s+1},\tau_t-k \sum^{\ell+s}_{i=\ell+1} (i,t)\right]
= 0,
$$
equivalently,
$$
[\tilde{\psi}_{s+1},B_t] = \frac{k}{2} \sum^{\ell+s}_{i=\ell+1} [\tilde{\psi}_{s+1},(it)].
$$
This shows
$$
\left[\tilde{\psi}_{s+1} - \sum^{\ell+s}_{i=\ell+1} (i,t)
\tilde{\psi}_{s+1}(i,t),B_t\right] = 0.  \hspace*{20pt} \Box
$$

\begin{theorem}
Suppose $\alpha$ satisfies $(>,I)$, then
$$
\psi_I|X = \Pi\{\kappa_i(\alpha)-\kappa_j(\alpha): \ \ell+1 \leq i < j
\leq \ell+m\}PH.
$$
\end{theorem}

{\em Proof.} By 4.9, $PH \zeta_\alpha = \sigma_I \zeta_{\sigma_{I}\alpha}$;
thus it suffices to compute the coefficient of $\zeta_{\sigma_{I}\alpha}$
in $\sigma_I \psi_I \zeta_\alpha$.  We use the inductive framework from
Definition 4.10 and Theorem 4.11.  For fixed $s < m$, let $\sigma_0,
\sigma_1$ be the reversing permutations for the intervals $[\ell+1,
\ell+s], \ [\ell+1, \ell+s+1]$ respectively. The inductive hypothesis is
that
$$
\psi_s \zeta_\gamma = \prod_{\ell+1 \leq i < j \leq \ell+s}
(\kappa_i(\gamma)-\kappa_j(\gamma))\sigma_0 \zeta_{\sigma_{0}\gamma}
$$
for any $\gamma$ satisfying $(>, [\ell+1, \ell+s])$ (trivial for $s =
1$).  Fix $\alpha$ satisfying $(>,I)$; let
$$
\pi_r := \prod \{\kappa_i(\alpha)-\kappa_j(\alpha): \ \ell+1 \leq i < j
\leq \ell+s +1, \ i \neq r, \ j \neq r\},
$$
and
$$
\pi_r' := \prod\{\kappa_i(\alpha): \ \ell+1 \leq i \leq \ell+s+1, \ i
\neq r\}.
$$
For $\ell+1 \leq j \leq \ell+s+1$, let $w_{(j)}$ be the cycle
$(\ell+s+1, \ell+s, \ldots,j+1,j)$ in $S_{[\ell+1,\ell+s+1]}$.  In the
notation of Definition 4.10,
$$
\psi_{s+1} = \sum^{\ell+s+1}_{j=\ell+1} (-1)^{\ell+s+1-j} w^{-1}_{(j)} \tilde{\psi}_{s+1}w_{(j)},
$$
because $w_{(j)} = (\ell+s,\ell+s-1,\ldots,j)(j,\ell+1+s)$ (in cycle
notation) and $\mbox{sgn}(w_{(j)}) = \ell+s+1-j$; coming from the skew
property of $\tilde{\psi}_{s+1}$ for $S_{[\ell+1,\ell+s]}$.

Let $\alpha_{(j)} = w_{(j)}\alpha$; thus
$$
\alpha_{(j)} = (\alpha_1,\ldots,\alpha_{\ell+1},\alpha_{\ell+2},\ldots,
\alpha_{j-1}, \alpha_{j+1}, \ldots,\alpha_{\ell+s+1}, \alpha_j,
\alpha_{\ell+s+2}, \ldots),
$$
which satisfies $(>, [\ell+1,\ell+s])$.  We claim the coefficient of
$\zeta_{\sigma_{1}\alpha}$ in $\sigma_1 w_{(j)}^{-1} \tilde{\psi}_{s+1}
w_{(j)}\zeta_\alpha$ is $\pi_j \pi'_j$.  By a standard identity for
alternating polynomials
$$
\sum^{\ell+s+1}_{j=\ell+1} (-1)^{\ell+s+1-j} \pi_j \pi'_j =
\prod_{\ell+1 \leq i < j \leq \ell+s+1} (\kappa_i(\alpha)-\kappa_j(\alpha)).
$$
To establish the claim:
$$
w_{(j)} \zeta_\alpha = \zeta_{\alpha_{(j)}} + \sum \left\{b(\gamma)
\zeta_\gamma: \ \gamma \in S_{[\ell+1,\ell+s+1]}\alpha \mbox{ and } \gamma
\succ \alpha_{(j)}\right\}
$$
(the triangularity of the matrix relating the bases $\{\zeta_{w\alpha}\}$
and $\{w \zeta_\alpha\}$).  By the inductive hypothesis and
$$
\sigma_1 w^{-1}_{(j)} \tilde{\psi}_{s+1} w_{(j)} \zeta_\alpha = \sigma_1
w_{(j)}^{-1} \sigma_0(\pi_j \pi'_j \zeta_{\sigma_{0} \alpha_{(j)}} +
\sum\{b'(\gamma) \zeta_{\sigma_{0}\gamma}: \ \gamma \succ
\alpha_{(j)}\}),
$$
for some coefficients $b'(\gamma)$.  But
$$
\sigma_1 w^{-1}_{(j)} \sigma_0
= w^{-1}_{(\sigma_{1}j)} \qquad (\sigma_1(j) = 2 \ell+s+2-j),
$$
which fixes $[\ell+1, 2\ell+s+1-j]$ pointwise, which shows that
$$
\sigma_1 w^{-1}_{(j)} \sigma_0 \zeta_{\sigma_{0}\gamma} \in
\mbox{span}\{\zeta_{w \sigma_{0}\gamma}: \ w \in S_{[2\ell+s+2-j,\ell+s+1]}\}.
$$
Now $(\sigma_0\gamma)_{\ell+s+1} = \gamma_{\ell+s+1} \in
\{\alpha_{j+1},\ldots,\alpha_{\ell+s+1}\}$, because $\gamma \succ
w_{(j)}\alpha$;  hence for any $\beta = w \sigma_0 \gamma$ in this span,
there must be an element of $\{\alpha_{j+1},\ldots,\alpha_{\ell+s+1}\}$
not appearing in $\{\beta_{\ell+1},\ldots,\beta_{2\ell+s+1-j}\}$, thus
$\beta \neq \sigma_1\alpha$.

The argument will be finished once it is shown that the coefficient of
$\zeta_{\sigma_{1}\alpha}$ in the $\{\zeta_{w\alpha}\}$ expansion of $w^{-1}
\zeta_{w\sigma_{1}\alpha}$ is 1: in the notation of (3.1),
\begin{eqnarray*}
w^{-1} \zeta_{w\sigma_{1}\alpha}
& = & w^{-1}\left(g_{w\sigma_{1}\alpha}+\sum_{\gamma \succ w
\sigma_{1}\alpha} B(\gamma, w \sigma_1 \alpha)g_\gamma\right) \\[4mm]
& = & g_{\sigma_{1}\alpha} + \sum_\gamma B(\gamma,w \sigma_1 \alpha) g_{w^{-1}\gamma};
\end{eqnarray*}
$\gamma \succ w \sigma_1 \alpha$ implies $w^{-1} \gamma \neq \sigma_1
\alpha$, the minimality of $\sigma_1 \alpha$ shows the coefficient of
$\zeta_{\sigma_{1}\alpha}$ in the right-hand side is 1. \hfill $\Box$
\bigskip

\begin{corollary}
$$
\psi_I a_{\alpha;I} = \prod_{\ell+1 \leq i < j \leq \ell+m}
(\kappa_i(\alpha) - \kappa_j(\alpha)-k)j_{\alpha,I}.
$$
\end{corollary}

{\em Proof.}  $\prod_{i < j}
(\kappa_i(\alpha)-\kappa_j(\alpha))\calE_-(\alpha^R;I)$ has the
specified value; see Lemma 3.13.  \hfill $\Box$
\bigskip

We return to the problem of evaluating $\calA_I f(1^N)$.  The proof of
the following comes later.

\begin{theorem}
Suppose $\{I_1,I_2,\ldots,I_t\}$ is a collection of disjoint
subintervals of $[1,N]$, $m_i := \# I_i$, $1 \leq i \leq t$ and $f$ is a
polynomial, then
$$
(\calA_{I_{1}} \calA_{I_{2}}\cdots A_{I_{t}}f)(1^N) =
\frac{1}{(Nk+1)_\mu \prod^t_{i=1} m_i!} \ (\psi_{I_{1}} \psi_{I_{2}}
\cdots \psi_{I_{t}} f)(1^N),
$$
where $\mu = (m_1-1,m_1-2,\ldots,1,0,m_2-1,\ldots,0,\ldots,m_t-1,\ldots,1,0,\ldots)^+.$
\end{theorem}

In the sequel, for an operator $A$ on polynomials, $A^*$ denotes the
adjoint with respect to the $p$-inner product.

Let $\iota(x) := \prod^N_{i=1} (1-x_i)^{-(Nk+1)}$, then $\langle
f,\iota\rangle_p = f(1^N)$ for any polynomial $f$; more generally, $\langle
f, F_k(\cd,z)\rangle_p = f(z)$.  Let $u_i := x_i/(1-x_i), \ 1 \leq i
\leq N$.

\begin{lemma}
$$
\calA^*_{I_{1}} \cdots \calA^*_{I_{t}} \iota(x) =
\prod^t_{i=1} \left(\frac{1}{m_i!} \ a_{I_{i}}(u)\right) \prod^N_{j=1} (1-x_j)^{-(Nk+1)}.
$$
\end{lemma}

{\em Proof.} First apply $\calA_{I_{j}}$ to $F_k(x,z)$ with respect to
$z$.  Without loss of generality, assume $I_j = [1,m]$ (with $m = m_j$); then
\begin{eqnarray*}
\calA_{I_{j}} F_k(x,z)
& = & \frac{1}{m!}\left(\sum_{w \in S_{[1,m]}}
\frac{\mbox{sgn}(w)}{\prod^m_{i=1} (1-x_i(zw)_i)}\right)
\prod^N_{j=m+1} (1-x_jz_j)^{-1}
\frac{\prod^N_{i,j=1} (1-x_iz_j)^{-k}}{a_{[1,m]}(z)} \\[4mm]
& = & \frac{1}{m!} a_{[1,m]} (x) \prod^m_{i,j=1} (1-x_iz_j)^{-1}
\prod^N_{j=m+1}(1-x_jz_j)^{-1} \prod^N_{i,j=1} (1-x_iz_j)^{-k},
\end{eqnarray*}
where the sum evaluates to $a_{[1,m]}(x)
a_{[1,m]}(z)/\prod^m_{i,j=1}(1-x_iz_j)$ by an identity of Cauchy (this
identity was used by Baker and Forrester [BF4] in a similar
calculation).  Because the intervals are disjoint,
$$
\calA_{I_{1}} \calA_{I_{2}} \cdots \calA_{I_{t}} F_k(x,z) =
\prod^t_{r=1} h_r(x,z) F_k(x,z),
$$
where
$$
h_r(x,z) = \frac{1}{m_r!} \ a_{I_{r}}(x) \prod\{(1-x_iz_j)^{-1}: \ i,j
\in I_r, \ i \neq j\}.
$$
To find $\calA_{I_{1}}^* \cdots \calA_{I_{t}}^* \iota$ put $z = 1^N$ in this
formula; note
$$
h_r(x,1^N) = \frac{1}{m_r!} a_{I_{r}}(x) \prod_{i \in I_{r}}
(1-x_i)^{-(m_i-1)} = \frac{1}{m_r!} a_{I_{r}}(u). \hspace*{20pt}
$$
\bigskip

The operators $\psi_I$ are generated by $\{T_i \rho_i: \ i \in I\}$ and
transpositions.  Recall $T_i \rho_i = T_i x_i +k$.  Write $T^u_i$ to
denote action in the variable $(u_1, u_2,\ldots)$.

\begin{lemma}
Suppose $f$ is a polynomial in $u$, then
$$
(T_i x_i+k)(f(u)\iota(x)) = ((1+u_i)(T^u_i u_i +k)f(u)) \iota(x), \ 1 \leq i
\leq N.
$$
\end{lemma}

{\em Proof.}  The product rule $T_i(h(x)\iota(x)) = (T_ih(x))\iota(x)+h(x)
\frac{\partial \iota(x)}{\partial x_i}$ applies because $\iota(x)$ is $S_N$-invariant.
The chain rule implies
\begin{eqnarray*}
\frac{[(T_ix_i+k)(f(u)\iota(x))]}{\iota(x)}
& = & f(u)(1+k)+(Nk+1)u_if(u)+u_i(1+u_i) \frac{\partial f(u)}{\partial
u_i} \\[4mm]
&& + \ k \sum_{j \neq i} \frac{u_i(1+u_j)f(u)-u_j(1+u_i)f(u(ij))}{u_i(1+u_j)-u_j(1+u_i)}
\end{eqnarray*}
(note $x_i = u_i/(1+u_i)$).  The typical term in the sum equals
$$
(1+u_i) \left(\frac{u_if(u)-u_jf(u(ij))}{u_i-u_j}\right) -u_i f(u).
\hspace*{20pt} \Box
$$
\bigskip

The effect on $f$ is to raise the degree by 1; in fact, the highest
degree term of the operator coincides with $T_i^{u*} = u_i(T_i^u u_i+k)$
(\cite{D3}, Proposition 4.1).

\begin{lemma} For an interval $I = [\ell+1,\ell+m]$ and a polynomial
$f(u)$ of degree $t$,
$$
\psi^*_I(f(u)\iota(x)) = \left[\prod_{\ell+1 \leq i < j \leq \ell+m}
(T^*_i-T^*_j)f(u)+f_1(u)\right]\iota(x),
$$
where $f_1(u)$ is a polynomial of degree $< t +m(m-1)/2$.
\end{lemma}

{\em Proof.}
Using the inductive framework from Theorem 4.11, suppose the statement
is true for the interval $[\ell+1,\ell+s]$ (trivial for $s = 1$). In the
notation of Definition 4.10, $\tilde{\psi}^*_{s+1} =
\psi^*_s(U_{\ell+1}\cdots U_{\ell+s})^*$ and by Lemma 4.16,
$$
U^*_i(f(u)\iota(x)) = ((1+u_i)(T_i^uu_i+k)f(u)-k \sum_{j < i} f(u(ij)))\iota(x),
$$
which has the same highest degree terms as $(T^*_if(u))\iota(x)$. By the
commutativity of $\{T^*_i\}$, the highest degree term of
$\tilde{\psi}^*_{s+1}(f(u)\iota(x))$ is
$$
\left(\prod^{\ell+s}_{i=\ell+1} T^*_i \prod_{\ell+1 \leq i < j \leq
\ell+s} (T^*_i-T^*_j)f(u)\right)\iota(x)
$$
(inductive hypothesis). Finally
$$
\psi^*_{s+1} = \tilde{\psi}^*_{s+1}-\sum^{\ell+s}_{i=\ell+1}
(i,\ell+s+1) \tilde{\psi}^*_{s+1}(i,\ell+s+1);
$$
and the usual identity for $a_{[\ell+1,\ell+s+1]}$ finishes the proof;
since
$$
(i,\ell+s+1)T^*_i(i,\ell+s+1) = T^*_{\ell+s+1}. \hspace*{20pt} \Box
$$

\begin{lemma}
$$
\left(\prod^t_{i=1} \psi_{I_{i}}^*\right) \iota(x) = c_I \prod^t_{i=1} a_{I_{i}}(u)
\iota(x)
$$
for some constant $c_I$.
\end{lemma}

{\em Proof.} By the previous lemma,
$$
\frac{\left(\prod^t_{i=1} \psi^*_{I_{i}}\right)\iota(x)}{\iota(x)} = f_0(u)+f_1(u),
$$
where $f_0(u) = \prod^t_{i=1} a_{I_{i}}(T^*)1$ and deg $f_1 <
\sum^t_{i=1} m_i(m_i-1)/2$. Also $f_0(u) +f_1(u)$ is skew for
$\prod^t_{i=1} S_{I_{i}}$ (direct product) by the skew property of
$\psi_I$; which implies $f_1 = 0$ and $f_0(u) = c_I\prod^t_{i=1}
a_{I_{i}}(u)$ because this is the unique skew polynomial of minimum
degree. \hfill $\Box$
\bigskip

{\em Proof of Theorem 4.14.} The lemmas show that
$$
\prod^t_{i=1} \psi^*_{I_{i}} \iota(x) =
\left(\left(\prod^t_{i=1} a_{I_{i}}(T^{u*})\right)1\right)\iota(x).
$$
In (\cite{D3}, Theorem 3.1) it was shown that $\prod_i a_{I_{i}} \in
E_\mu$ for $\mu = (m_1-1, m_1-2, \ldots,1,0,m_2-1,\ldots,1,0,\ldots)^+$,
thus (by Corollary 2.6)
$$
\prod^t_{i=1} a_{I_{i}}(T^{u*})1 = (Nk+1)_\mu \prod^t_{i=1} a_{I_{i}}(u)
$$
(see also \cite{DH}); and so
$$
\prod^t_{i=1} \psi^*_{I_{i}}\iota = (Nk+1)_\mu\left(\prod^t_{i=1} m_i!\right)
\prod^t_{i=1} \calA^*_{I_{i}} \iota. \hspace*{20pt}
$$

\begin{corollary}
Suppose $\alpha$ satisfies $(>,I_i)$ for each $i, \ 1 \leq i \leq t$, then
\begin{eqnarray*}
\lefteqn{\hspace*{-30pt}
\left(\prod^t_{i=1} \calA_{I_{i}} \sum_{w \in S_{I}} (\mbox{sgn }w) w
\zeta_\alpha\right) (1^N)} \\[4mm]
& = & \frac{(Nk+1)_{\alpha^+}}{(Nk+1)_\mu} \
\frac{\calE_-(\alpha,[1,N])}{h(\alpha^+,1)} \ \prod^t_{r=1} \
\prod\{\kappa_i(\alpha)-\kappa_j(\alpha)-k: \ 
i,j \in I_r, \ i < j\},
\end{eqnarray*}
where $S_I := S_{I_{1}} \times S_{I_{2}} \cdots \times S_{I_{t}}$.
\end{corollary}

{\em Proof.} Let $\pi_r = \prod\{\kappa_i(\alpha)-\kappa_j(\alpha)-k: \
i,j \in I_r, \ i < j\}$.  For each interval $I_r$, $\psi_{I_{r}} \sum_{w
\in S_{I_{r}}} \mbox{sgn}(w) w \zeta_\alpha = \pi_r \sum_{w \in S_{I_{r}}}w
\zeta_\alpha$, by Corollary 4.13.  Since $[\psi_{I_{r}},w] = 0$ for $w \in
S_{I_{i}}$ and $[\psi_{I_{r}},\psi_{I_{i}}] = 0$, for $i \neq r$
(Theorem 4.11),
$$
\prod^t_{r=1} \psi_{I_{r}} \sum_{w \in S_{I_{r}}} \mbox{sgn}(w) w
\zeta_\alpha = \prod^t_{r=1} \pi_r \sum_{w \in S_{I}} w \zeta_\alpha
$$
(the sum is a product of sums over $S_{I_{1}},\ldots,S_{I_{t}}$).  Also
$$
\sum_{w \in S_{I}} w \zeta_\alpha (1^N) = \left(\prod^t_{i=1} m_i!\right)
\zeta_\alpha (1^N),
$$
and
$$
\zeta_\alpha (1^N) = (Nk+1)_\alpha \calE_-(\alpha,[1,N])/h(\alpha^+,1)
$$
(by (3.2)). Further
$$
\calA_{I_{r}} \sum_{w\in S_{I_{r}}} \mbox{sgn}(w) w\zeta_\alpha =
\sum_{w\in S_{I_{r}}} \mbox{sgn}(w) w \zeta_\alpha/a_{I_{r}}
$$
(and this is a typical term in the product sum over $S_{I_{1}} \times
\cdots \times S_{I_{t}}$). \hfill $\Box$

\section{Orthogonal Polynomials of Type $B_N$}
\setcounter{equation}{0}

This section deals with operators and orthogonal decompositions
associated with the group generated by sign-changes and permutations of
coordinates. Previously, Baker and Forrester \cite{BF2} considered some
of the orthogonal polynomials, called generalized Laguerre polynomials.
These come from polynomials which are even in each coordinate or odd in
each coordinate. The general situation is developed in the sequel. One
consequence is a complete set of eigenfunctions for the Hamiltonian of
the $B_N$ spin Calogero model ($1/r^2$ interactions confined in harmonic
potential), with arbitrary parity, that is, for any subset $A \subset
[1,N]$ we find the eigenfunctions which are odd in $x_i$, $i \in A$ and
even in $x_i$, $i \notin A$.

The underlying symmetry group $W_N$ is the Weyl group of type $B_N$,
called the hyperoctahedral group.  It is generated by permutations of
coordinates and sign-changes on $\mathbb{R}^N$.  The reflections in $W_N$
are $\{\sigma_{ij}, \tau_{ij}: \ 1 \leq i < j \leq N\}$ and $\{\sigma_i:
\ 1 \leq i \leq N\}$, defined by $x \sigma_{ij} = (x_1,\ldots
\stackrel{i}{x}_{j},\ldots,\stackrel{j}{x}_i,\ldots)$, $x \tau_{ij} =
(x_1, \ldots, -\stackrel{i}{x}_j,\ldots, -\stackrel{j}{x}_i,\ldots)$, $x
\sigma_i = (x_1, \ldots, -\stackrel{i}{x}_i,\ldots)$.  There are two
parameters $k$, $k_1$ in the algebra of differential-difference operators:
\begin{eqnarray}
T^B_i f(x)
& := & \frac{\partial f}{\partial x_i} + k_1
\frac{f(x)-f(x\sigma_i)}{x_i} \\[4mm]
&& + \ k \sum_{j \neq i} \left\{\frac{f(x)-f(x\sigma_{ij})}{x_i-x_j} + \frac{f(x)-f(x\tau_{ij})}{x_i+x_j}\right\}
\nonumber
\end{eqnarray}
($1 \leq i \leq N$, for convenience $\sigma_{ij} = \sigma_{ji}$,
$\tau_{ij} = \tau_{ji}$ for $j < i$).

The $S_N$-theory can be applied to the analysis of $W_N$ by writing
polynomials in the form $x_A g(x^2_1, x^2_2, \ldots,x^2_N)$ where $x_A
:= \prod_{i\in A} x_i$, $A \subset [1,N]$.  For a composition $\alpha$
let $\hat{p}_\alpha (x) := x_A p_\beta(x^2_1,\ldots,x^2_N)$ where $A =
\{i: \ \alpha_i \mbox{ is odd}\}$, and $\beta_i = \lfloor
\alpha_i/2\rfloor$, $1 \leq i \leq N$.  Observe
$$
\sum_\alpha
\hat{p}_\alpha(x) z^\alpha = \prod^N_{i=1} ((1-x_iz_i)^{-1}
\prod^N_{j=1} (1-x^2_i z_j^2)^{-k}).
$$
The raising operator $\hat{\rho}_i$ is defined by $\hat{\rho}_i
\hat{p}_\alpha = \hat{p}(\alpha_1,\ldots,\alpha_{i+1},\ldots)$.

We will use $y \in \mathbb{R}^N$ to denote $(x^2_1,x^2_2,\ldots,x^2_N)$ and
the $A$-type operators $T_i$ act on $y$.  The following properties hold
(\cite{D4}, Prop. 2.1), for $A \subset [1,N]$, any polynomial $g$ in $y$:
\begin{eqnarray}
(\sigma_{ij}+\tau_{ij}) x_A g(y)
& \!\! = \!\! & \left\{\begin{array}{ll}
(2x_A)(ij)g(y) & i,j \in A \mbox{ or } i,j \notin A, \\[2mm]
0 & \mbox{else};
\end{array}\right. \\[5mm]
T^B_i \hat{\rho}_i(x_A g(y))
& \!\! = \!\! & \left\{\begin{array}{ll}
2 x_A(T_i \rho_i g(y)), & i \in A, \\[2mm]
2x_A((k_1-k-\frac{1}{2})g(y)+T_i \rho_i g(y)-k \sum_{s\in A}(ij)g(y)),
& i \notin A.
\end{array}\right.
\end{eqnarray}

\begin{definition}
For $1 \leq i \leq N$, $U^B_i := T^B_i \hat{\rho}_i-k \sum_{j
<i}(\sigma_{ij}+\tau_{ij})$ acts on polynomials in $x$. For a subset $A
\subset [1,N]$, let
$$
U_{A,i} := T_i \rho_i-k \sum\{(ij): \ j < i \mbox{
and } j \in A\} \mbox{ for } i \in A,
$$
$$
U_{A,i} := T_i
\rho_i+(k_1-k-\frac{1}{2})1-k \sum\{(ij): \ j \in A
\mbox{ or } (j \notin A \mbox{ and } j < i)\} \mbox{ for }
i \notin A,
$$
acting on polynomials in $y$.
\end{definition}

\begin{proposition}
For $A \subset [1,N]$, $U^B_i (x_A g(y)) = 2 x_A U_{A,i}g(y), \ 1 \leq i
\leq N$. Also
$$
\prod_{i \in A} T^B_i (x_A g(y)) = 2^{\#A} \prod_{i \in A} \left(U_{A,i}+k_1-k-\frac{1}{2}\right)g(y).
$$
\end{proposition}

This proves the commutativity of $\{U^B_i\}$ in terms of propositions
about $S_N$; a direct proof is also possible.

Let $V^B$ denote the intertwining operator for $W_N$, thus $T^B_i V^B =
V^B \frac{\partial}{\partial x_i}$, and $V^B$ is homogeneous of degree
0; also let
$\xi^B$ be the linear operator on polynomials defined by
$\xi^B \hat{p}_\alpha = x^\alpha/\alpha!$. Here is a description of
the eigenspace decomposition of $V^B
\xi^B$ constructed in \cite{D4}.
Each space is irreducible for the algebra generated by $\{T^B_i
\hat{\rho}_i: \ 1 \leq i \leq N\}$.

\begin{definition}
A $B$-partition $\alpha$ is a composition $(\alpha_1, \alpha_2,
\ldots,\alpha_N)$ where the odd and even parts are respectively
nonincreasing, that is, $i < j$ and $\alpha_i \equiv \alpha_j$ {\rm mod} 2
implies $\alpha_i \geq \alpha_j$.  A standard $B$-partition is one in
which the odd parts come first (for some $\ell$, $\alpha_i$ is odd for $i
\leq \ell$, even for $i > \ell$). For $\alpha \in \calN_N$, let
$$
h(\alpha) = (\lfloor \alpha_1/2\rfloor, \lfloor
\alpha_2/2\rfloor,\ldots), \ b(\alpha) = (\alpha_1 -\lfloor
\alpha_1/2\rfloor, \alpha_2-\lfloor\alpha_2/2\rfloor, \ldots).
$$
\end{definition}

This is an example of a $B$-partition: $\alpha = (4,5,3,4,2,0,1)$, and
$h(\alpha) = (2,2,1,2,1,0,0)$, $b(\alpha) = (2,3,2,2,1,0,1)$.  The
corresponding standard $B$-partition is $(5,3,1,4,4,2,0)$.

Any $B$-partition can be rearranged to a standard one; if $\alpha$ is a
$B$-partition with $\ell$ odd parts, then there exists a standard
$B$-partition $\tilde{\alpha}$ and $w \in S_N$ so that $w\tilde{\alpha}
= \alpha$ and $1 \leq i < j \leq \ell$ or $\ell+1 \leq i < j \leq N$
implies $w(i) < w(j)$ (note $\alpha_{w(i)} = \tilde{\alpha}_i$, each $i$).

The property of permutations described above will be used several times
in the sequel.

\begin{definition}
For any subset $A \subset [1,N]$ let $w_A \in S_N$ be the unique
permutation satisfying: $w_A([1,\ell]) = A \quad (\mbox{for } \ell =
\#A),$ and $1 \leq i < j \leq \ell$ or $\ell+1 \leq i < j \leq N$
implies $w_A(i) < w_A(j)$. If $1 \leq \ell < N$, the correspondence $A
\rarrow w_A$ is one-to-one.
\end{definition}

\begin{proposition}
Let $A \subset [1,N]$, $\ell = \#A$, then
\begin{eqnarray*}
U^B_i w_A(x_1x_2\cdots x_\ell g(y))
& = & w_A(2 x_1 \cdots x_\ell U_{[1,\ell],s} g(y)) \\[4mm]
& = & w_A(U^B_s (x_1 \cdots x_\ell g(y))), \mbox{ for } 1 \leq i \leq N,
\ s = w^{-1}(i).
\end{eqnarray*}
\end{proposition}

{\em Proof.} $(T^B_i \hat{\rho}_i)w_A = w_A T^B_s \hat{\rho}_s$, for $s
= w_A^{-1}(i)$ and $(ij)w_A = w_A(w_A^{-1}(i),w_A^{-1}(j))$ for $j \neq
i$. The order preserving properties of $w_A$ and Proposition 5.5 imply
the stated equations. \hfill $\Box$
\bigskip

The proposition shows how to find joint eigenfunctions of $\{U_i^B\}$
corresponding to arbitrary $B$-partitions once the standard case is done.

For a given standard $B$-partition $\alpha$, with $\alpha_i$ being odd
exactly when $1 \leq i \leq \ell$, let $G_\ell = S_{[1,\ell]} \times
S_{[\ell+1,N]}$, $\beta = h(\alpha), \ \mu = \beta^+$.  Define
$$
E^B_\alpha = \mbox{span}\{x_1x_2 \cdots x_\ell \zeta_\gamma(y): \, \gamma
= w \beta, \ w \in G_\ell\}.
$$

\begin{proposition} $E^B_\alpha$ is invariant under $G_\ell$, $U^B_i
(x_1 \cdots x_\ell \zeta_\gamma(y)) = c_i x_1 \cdots x_\ell
\zeta_\gamma(y)$, where $c_i = 2 \kappa_i (\gamma)$ for $i \in [1,\ell]$,
$c_i = 2(\kappa_i(\gamma)+k_1-k-\frac{1}{2})$ for $i \in [\ell+1,N]$.
\end{proposition}

This follows from Proposition 5.2. To express the eigenvalues of $V^B
\xi^B$, let
$$
\Lambda(\alpha) := (Nk+1)_{h(\alpha)^+}((N-1)k+k_1+\frac{1}{2})_{b(\alpha)^+},
$$
for $ \alpha \in \calN_N$.  Then $E^B_\alpha$ is an eigenspace for $V^B
\xi^B$ with eigenvalue $2^{-|\alpha|} \Lambda(\alpha)^{-1}$.  It was
shown in (\cite{D4}, Proposition 4.2) that $x_1 \cdots x_\ell
\zeta_\beta(y)$ is an eigenfunction of $V^B
\xi^B$ with this eigenvalue,
and $V^B\xi^B$ commutes with $w \in W_N$.

Again there are three inner products in which each $U^B_i$ is
self-adjoint, and which are $W_N$-invariant: (1) the $p$-product,
$\langle x^\alpha, \hat{p}_\beta\rangle_p = \delta_{\alpha\beta}$; (2)
the $B$-product, $\langle f,g\rangle_B = f(T^B)g(x)|_{x=0}$ (depends on
both parameters); (3) the $N$-torus product,
$$
\langle f,g\rangle_{\mathbb{T}} := c_k \int_{\mathbb{T}^N} f(x)
\check{g}(x) \prod_{1 \leq j < \ell \leq
N}|(x^2_j-x^2_\ell)(x^{-2}_j-x_\ell^{-2})|^k dm(x)
$$
(same notation as in 3.20; $c_k$ is the normalizing constant).  It is obvious
that $\langle x_{A_{1}}g_1(y), x_{A_{2}}g_2(y)\rangle = 0$ if $A_1, A_2
\subset [1,N]$ and $A_1 \neq A_2$, by the $W_N$-invariance (change the
sign of $x_i$ for $i \in A_1\setminus A_2$ or $i \in A_2 \setminus
A_1$).  Further $\langle x_A g_1(y), x_Ag_2(y)\rangle = \langle g_1,
g_2\rangle$ for the $p$- or $\mathbb{T}$-products, using the type A
definition on the right.

By the results in Section 3, let $f(x) = x_1 \cdots x_\ell \zeta_\gamma(y)
\in E^B_\alpha$, then
$$
\|f\|^2_p = \|\zeta_\gamma\|^2_p = \calE_+(\gamma)\calE_-(\gamma)h(\gamma^+,k+1)/h(\gamma^+,1).
$$
By the same argument as Proposition 2.5,
\beq
\|f\|^2_B = f(T^B)f(x) = 2^{|\alpha|}\Lambda(\alpha)\|\zeta_\gamma\|^2_p.
\eeq

In \cite{D2} we showed that
\beq
\langle f,g\rangle_B = c_{k,k'} \int_{\mathbb{R}^N} \left(e^{-L/2}f\right)
\left(e^{-L/2}g\right) \prod^N_{j=1} |x_j|^{2k_1} \prod_{1 \leq i < j
\leq N} |x^2_i-x^2_j|^{2k} e^{-|x|^{2}/2}dx,
\eeq
for polynomials $f, g$ where $L := \sum^N_{i=1}(T^B_i)^2$ (the
normalizing constant $c_{k,k'}$ is chosen to make $\langle 1,1\rangle_B
= 1$ and is computed by means of the Macdonald-Selberg integral). 
Define a generalized Hermite polynomial
\beq
H_\beta(x) := e^{-L/2}(x_1\cdots x_\ell \zeta_\gamma(y)),
\eeq
with
$$
\beta = (2\gamma_1+1,\ldots,2\gamma_\ell+1,
2\gamma_{\ell+1},\ldots,2\gamma_N) \quad
(\mbox{so } \gamma = h(\beta)).
$$
The complete orthogonal basis is given by $\{w_A H_\beta(x): \ A \subset
[1,N], \ \ell = \#A, \ \beta_i$ is odd exactly for $i \in [1,\ell]\}$.

By using the ``non-symmetric'' binomial coefficients introduced by Baker
and Forrester \cite{BF3} we can produce an expression for $H_\beta$ in
terms of $\{\zeta_\gamma\}$; although very little is known about the coefficients.

\begin{definition}
For $\alpha \in \calN_N$, the binomial coefficients (depending on $k$)
are implicitly defined by
$$
\frac{\zeta_\alpha(y+1^N)}{\zeta_\alpha(1^N)} = \sum_\gamma \left(\alpha
\atop \gamma\right) \frac{\zeta_\gamma(y)}{\zeta_\gamma(1^N)}, \ y \in \mathbb{R}^N.
$$
\end{definition}

It is known that $\left(\alpha \atop \gamma\right) = 0$ unless $\gamma^+
\subset \alpha^+$ (that is, $(\gamma^+)_i \leq (\alpha^+)_i, \ 1 \leq i
\leq N$).  For any scalar $s$, the homogeneity of $\zeta_\alpha$ implies
$$
\frac{\zeta_\alpha(y+s 1^N)}{\zeta_\alpha (1^N)} = \sum_\gamma \left(\alpha
\atop \gamma\right) s^{|\alpha|-|\gamma|}
\frac{\zeta_\gamma(y)}{\zeta_\gamma(1^N)}, \ y \in \mathbb{R}^N.
$$
When $k = 0$, $\left(\alpha \atop \gamma\right) = \prod^N_{i=1}
\left(\alpha_i \atop \gamma_i\right)$ (ordinary binomial coefficients).

\begin{proposition}
(Baker and Forrester \cite{BF3}) For $\alpha \in \calN_N$,
$$
\exp\left(s \sum^N_{i=1} y_i\right) \zeta_\alpha(y)
= \sum_{\gamma^+ \supset \alpha^+}
\frac{h(\alpha^+,k+1)\calE_+(\alpha)}{h(\gamma^+,k+1)\calE_+(\gamma)}
\left(\gamma \atop \alpha\right) s^{|\gamma|-|\alpha|} \zeta_\gamma (y), \
y \in \mathbb{R}^N, \ s \in \mathbb{R}.
$$
\end{proposition}

{\em Proof.} The adjoint in the $A$-inner product of multiplication by
$\exp(s \sum_i y_i)$ is translation $g(y) \mapsto g(y+s 1^N)$. Indeed,
suppose $f$ and $g$ are polynomials, then
\begin{eqnarray*}
\langle
e^{s \Sigma y_i} f(y),g(y)\rangle_A
& = & f(T) e^{s \Sigma T_i} g(y)|_{y=0} \\[3mm]
& = & f(T) \exp\left(s \sum^N_{i=1} \frac{\partial}{\partial y_i}\right)
g(y)|_{y=0} \\[3mm]
& = & f(T) g(y+s 1^N)|_{y=0},
\end{eqnarray*}
because $\sum^N_{i=1} T_i = \sum^N_{i=1} \frac{\partial}{\partial y_i}$.
 The given expression is found by using the $A$-norms of the orthogonal
basis elements $\zeta_\gamma$ (see Corollary 3.19). \hfill $\Box$
\bigskip

The adjoint of $e^{-sL}$ in the $B$-product is multiplication by
$\exp\left(-s \sum^N_{i=1} x^2_i\right)$, which can be evaluated using
the previous result for $\zeta_\gamma(y)$. For $\beta \in \calN_N$, $0
\leq \ell \leq N$, let
$$
b(\beta,\ell) = (2 \beta_1+1,\ldots,2 \beta_\ell+1, 2\beta_{\ell+1}, \ldots,2\beta_N).
$$
For $\alpha \in \calN_N, \ s \in \mathbb{R}$,
\beq
e^{sL}(x_1\cdots x_\ell \zeta_\alpha(y)) = \sum_{\beta^+ \subset \alpha^+}
\frac{\Lambda(b(\alpha,\ell))
h(\beta^+,1)\calE_-(\alpha)}{\Lambda(b(\beta,\ell))h(\alpha^+,1)\calE_-(\beta)}
\, (4s)^{|\alpha|-|\beta|} x_1 x_2 \cdots x_\ell \zeta_\beta(y).
\eeq

The formula follows from a similar adjoint-type calculation as in 5.8.
Here the norm of the orthogonal basis element is
\begin{eqnarray*}
\|x_1x_2\cdots x_\ell \zeta_\beta(y)\|^2_B
& = & 2^{2|\beta|+\ell} \Lambda(b(\beta,\ell))\|\zeta_\beta\|^2_p \\[4mm]
& = & 2^{2|\beta|+\ell} \Lambda(b(\beta,\ell))h(\beta^+,k+1)\calE_+(\beta)\calE_-(\beta)/h(\beta^+,1)
\end{eqnarray*}
(see (5.4)).  Put $s = - \frac{1}{2}$ in the formula to produce an
orthogonal basis for $e^{-|x|^{2}/2}$, $s = - \frac{1}{4}$ for
$e^{-|x|^{2}}$ in the formula (5.5); that is, including all the
polynomials $w_A(x_1\cdots x_\ell \zeta_\alpha(y)), \ A \subset [1,N], \
\ell = \#A$.

The special cases $\ell = 0$ and $\ell = N$ have already been obtained
by Baker and Forrester \cite{BF3}, who called them generalized Laguerre polynomials.

We discuss the connection to the $B_N$-type spin Calogero model in (1.3).

>From the commutation $[L,x_i] = 2T^B_i$ (\cite{D1}, Proposition 2.2), it
follows that $e^{-L/2} x_i = (x_i-T^B_i) e^{-L/2}$.  Also
$$
T^B_i x_i = x_i T_i^B+1+2k_1\sigma_i+k \sum_{j \neq i} (\sigma_{ij}+\tau_{ij}),
$$
which shows that
$$
e^{-L/2} U^B_i e^{L/2} = x_i T_i^B-(T_i^B)^2+1+k+(2k_1-k)\sigma_i +k
\sum_{j > i}(\sigma_{ij}+\tau_{ij}).
$$
The Hermite polynomials from (5.6) (with $s = -\frac{1}{2}$) are
simultaneous eigenfunctions of these operators.  Then
\begin{eqnarray*}
\calH_3 & := & e^{-L/2} \sum^N_{i=1} U^B_i e^{L/2} \\[4mm]
& = & \sum^N_{i=1} x_i \partial_i - \sum^N_{i=1} (T^B_i)^2+(k_1-k)
\sum^N_{i=1} \sigma_i +N(Nk+1+k_1).
\end{eqnarray*}
The eigenvalues depend only on the degree and the number of odd indices,
$$
\calH_3 e^{-L/2} (x_1 \cdots x_\ell \zeta_\alpha(x^2)) =
((2|\alpha|+\ell)+2\ell(k-k_1) +N(Nk+1+2k_1-k)(e^{-L/2} x_1 \cdots
x_\ell \zeta_\alpha(x^2)).
$$
Let
$$
h(x) = \prod_{1 \leq i < j \leq N} |x^2_i-x^2_j|^k \prod^N_{i=1} |x_i|^{k_1},
$$
then
\begin{eqnarray*}
h(x) \calH_3 h(x)^{-1}
& = & \sum^N_{i=1} \left(x_i \, \frac{\partial}{\partial x_i} -
\frac{\partial^2}{\partial x^2_i}\right) \\[4mm]
&& + \ (k_1-k) \sum^N_{i=1} \sigma_i + k_1 \sum^N_{i=1}
\frac{k_1-\sigma_i}{x^2_i} \\[4mm]
&& + \ 2k \sum_{1 \leq i < j \leq N}
\left\{\frac{k-\sigma_{ij}}{(x_i-x_j)^2} +
\frac{k-\tau_{ij}}{(x_i+x_j)^2}\right\} +N(k+1).
\end{eqnarray*}
Conjugating once more leads to
$$
e^{-|x|^{2}/4} h(x) \calH_3 h(x)^{-1} e^{|x|^{2}/4}
= \calH_2+(k_1-k) \sum^N_{i=1} \sigma_i +N(k+\frac{1}{2}).
$$
The middle term is basically counting the odd indices. Yamamoto \cite{Y} already
found that the eigenvalues were evenly spaced.  Baker and Forrester
\cite{BF3} studied the $A$-version of this model with a similar
transformation, namely $\exp \left(-\sum^N_{i=1}T^2_i/2\right)$. Van
Diejen \cite{vD} and Kakei \cite{K1}, \cite{K2}, \cite{K3} have studied
the symmetric $(W_N$-invariant) eigenfunctions of this model.

The theory of symmetric and alternating polynomials from Sections 3 and
4 can be applied. We will only write down the two-interval situation,
but the methods apply to finer partitions as well.

Fix an interval $[1,\ell]$, let $G_\ell = S_{[1,\ell]} \times
S_{[\ell+1,N]}$, and choose $\alpha \in \calN_N$ which satisfies $(\geq,
[1,\ell])$ and $(\geq, [\ell+1,N])$ (corresponding to a standard
$B$-partition). Then
$$
\alpha^R := \sigma_{[1,\ell]}\sigma_{[\ell+1,N]}\alpha = (\alpha_\ell,\alpha_{\ell-1},\ldots,\alpha_1,\alpha_N,\ldots,\alpha_{\ell+1}),
$$
and $\# G_\ell \alpha$ is the number of distinct permutations of
$(\alpha_1,\ldots,\alpha_\ell), (\alpha_{\ell+1},\ldots,\alpha_N)$.  The
following polynomial is an invariant of $S_{[1,\ell]} \times
W_{[\ell+1,N]}$: Let
\begin{eqnarray*}
j_{\alpha;\ell} & := & x_1 x_2 \cdots x_\ell \calE_+(\alpha^R,[1,\ell])
\calE_+(\alpha^R,[\ell+1,N]) \\[4mm]
&& \cdot \ \sum_{\beta \in G_{\ell}\alpha}
\frac{1}{\calE_+(\beta,[1,\ell])\calE_+(\beta,[\ell+1,N])} \
\zeta_\beta(x^2_1,\ldots, x^2_N).
\end{eqnarray*}
Then $j_{\alpha;\ell} = x_1x_2 \cdots x_\ell \sum_w w\zeta_\alpha$
(summing over a complete set of representatives for the cosets $w\{w_0
\in G_\ell: \ w_0 \alpha = \alpha\}$. Further,
$$
\|j_{\alpha;\ell}\|^2_p = (\# G_\ell
\alpha)\calE_+(\alpha^R,[1,\ell]) \calE_+(\alpha^R,[\ell+1,N])\|\zeta_\alpha\|^2_p,
$$
where
$$
\|\zeta_\alpha\|^2_p = \calE_+(\alpha)\calE_-(\alpha)h(\alpha^+,k+1)/h(\alpha^+,1);
$$
and
$$
\|j_{\alpha;\ell}\|^2_B = 2^{2|\alpha|+\ell} \Lambda(b(\alpha,\ell))\|j_{\alpha;\ell}\|^2_p.
$$
This is also the squared norm (from (5.5)) of $e^{-L/2}j_{\alpha;\ell}$,
which has the same $S_{[1,\ell]}\times W_{[\ell+1,N]}$ invariance. (For
an interval $I$, $W_I$ is the group generated by $S_I$ and $\{\sigma_i:
\ i \in I\}$.)

Further $j_{\alpha;\ell}(1^N) = (\#G_\ell \alpha)\calE_-(\alpha)(Nk+1)_{\alpha^+}
/h(\alpha^+,1)$.

Suppose that $\alpha$ satisfies $(>,[1,\ell])$ and $(>,[\ell+1,N])$. 
The following polynomial is alternating for $W_{[1,\ell]} \times
S_{[\ell+1,N]}$. Let
\begin{eqnarray*}
a_{\alpha;\ell}
& = & x_1 \cdots x_\ell \calE_-(\alpha^R,[1,\ell])
\calE_-(\alpha^R,[\ell+1,N]) \\[4mm]
&& \cdot \ \sum_{w \in G_{\ell}}
\frac{\mbox{sgn}(w)}{\calE_-(w\alpha,[1,\ell]) \calE_-(w
\alpha,[\ell+1,N])} \ \zeta_{w\alpha}(x^2_1,\ldots,x^2_N).
\end{eqnarray*}
Then $a_{\alpha;\ell} = x_1 \ldots x_\ell \sum_{w \in G_{\ell}}
\mbox{sgn}(w) w \zeta_\alpha$,
$$
\|a_{\alpha;\ell}\|^2_p = \ell! (N-\ell)! \calE_-(\alpha^R,[1,\ell])
\calE_-(\alpha^R,[\ell+1,N]) \|\zeta_\alpha\|^2_p,
$$
and
$$
\|a_{\alpha;\ell}\|^2_B = 2^{2|\alpha|+\ell} \Lambda(b(\alpha,\ell))\|a_{\alpha;\ell}\|^2_p.
$$
Further,
\begin{eqnarray*}
\lefteqn{\hspace*{-15pt}
a_{\alpha;\ell}(x)/(a_{[1,\ell]}(x^2)a_{[\ell+1,N]}(x^2))\bigg|_{x=1^N}} \\[4mm]
& = & \frac{(Nk+1)_{\alpha^+} \calE_-(\alpha)}{(Nk+1)_\mu h(\alpha^+,1)}
\prod \{\kappa_i(\alpha)-\kappa_j(\alpha)-k: \ 1 \leq i,j \leq \ell
\mbox{ or } \ell+1 \leq i < j \leq N\},
\end{eqnarray*}
where $\mu = (\ell-1, \ell-2, \ldots, 1, 0, N-\ell-1, \ldots, 1,0)^+$;
recall $a_{[1,\ell]}(x^2) = \prod_{1 \leq i < j \leq \ell} (x^2_i
-x^2_j)$. Again $e^{-L/2} a_{\alpha;\ell}$ is also alternating for
$W_{1,\ell]} \times S_{[\ell+1,N]}$, and its squared norm (for (5.5))
equals $\|a_{\alpha;\ell}\|^2_B$.

As mentioned above, the techniques of sections 2 and 3 can be used to
describe polynomials with prescribed symmetry for direct products
$W_{I_{1}} \times W_{I_{2}} \cdots \times W_{I_{r}} \times S_{I_{r+1}}
\cdots \times S_{I_{t}}$ for any collection $\{I_1, I_2, \ldots, I_t\}$
of pairwise disjoint intervals in $[1,N]$.

More information about the generalized binomial coefficients
needs to be obtained so that more concrete algorithms for the type-B
Hermite polynomials can be found. The polynomials could be useful in the
numerical cubature associated to the Macdonald-Mehta-Selberg integral.
Also, such knowledge could lead to orthogonal bases for harmonic
polynomials, which, by definition, are annihilated by
$\sum^N_{i=1}(T^B_i)^2$; for example, these polynomials appear when one
uses spherical polar coordinate systems to find eigenfunctions of some
Hamiltonians, see Section 3.4 in \cite{vD}.


\begin{thebibliography}{BF1}
\bibitem[BF1]{BF1}  Baker, T. H. and Forrester, P. J. : The
Calogero-Sutherland model and generalized classical polynomials; Commun. Math.
Phys. 188(1987), 195-216; preprint, solv-int/9608004.

\bibitem[BF2]{BF2}  Baker, T. H. and Forrester, P. J. : The
Calogero-Sutherland model and polynomials with prescribed symmetry; Nucl.
Phys. B492(1997), 682-716.

\bibitem[BF3]{BF3}  Baker, T. H. and Forrester, P. J. : Non-symmetric Jack
polynomials and integral kernels; Duke Math. J., to appear; preprint,
q-alg/9612003.

\bibitem[BF4]{BF4}  Baker, T. H. and Forrester, P. J. : Symmetric Jack
polynomials from the non-symmetric theory; preprint, q-alg/9707001, 1 Jul.
1997.

\bibitem[BDF]{BDF}  Baker, T. H., Dunkl, C. F. and Forrester, P. J. :
Polynomial eigenfunctions of the Calogero-Sutherland-Moser models with
exchange terms; Proc. CRM Workshop on Calogero-Sutherland-Moser models, to
appear.

\bibitem[BO]{BO}  Beerends, R. and Opdam, E. : Certain hypergeometric series
related to the root system $BC$; Trans. Amer. Math. Soc. 339(1993), 581-609.

\bibitem[C]{C}  Cherednik, I. : A unification of the Knizhnik-Zamolodchikov
and Dunkl operators via affine Hecke algebras; Inv. Math. 106(1991), 411-432.

\bibitem[vD]{vD}  van Diejen, J. F. : Confluent hypergeometric orthogonal
polynomials related to the rational quantum Calogero system with harmonic
confinement, Commun. Math. Phys. 188 (1997), 467--497.

\bibitem[D1]{D1}  Dunkl, C. F. : Differential-difference operators
associated to reflection groups, Trans. Amer. Math. Soc. 311 (1989), 167-183.

\bibitem[D2]{D2}  Dunkl, C. F. : Integral kernels with reflection group
invariance; Canadian J. Math. 43(1991), 1213-1227.

\bibitem[D3]{D3}  Dunkl, C. F. : Intertwining operators and polynomials
associated with the symmetric group; Monatsh. Math., to appear.

\bibitem[D4]{D4}  Dunkl, C. F. : Intertwining operators of type $B_N$; Proc.
CRM Workshop on $q$-special functions and algebraic methods; preprint
CRM-2380, 1996.

\bibitem[DH]{DH}  Dunkl, C. F. and Hanlon, P. : Integrals of polynomials
associated with tableaux and the Garsia-Haiman conjecture; Math. Z.228(1998),
537-567.

\bibitem[K1]{K1}  Kakei, S. : Common algebraic structure for the
Calogero-Sutherland models; preprint, solv-int/9608009, 27 Aug. 1996.

\bibitem[K2]{K2}  Kakei, S. : An orthogonal basis for the $B_N$ - type
Calogero model; preprint, solv-int/9610010, 28 Oct. 1996.

\bibitem[K3]{K3}  Kakei, S. : Intertwining operators for a degenerate double
affine Hecke algebra and multivariable orthogonal polynomials; preprint,
q-alg/9706019, 17 Jun. 1997.

\bibitem[LV1]{LV1}  Lapointe, L. and Vinet, L. : A Rodrigues formula for the
Jack polynomials and the Macdonald-Stanley conjecture; IMRN 9(1995), 419-424.

\bibitem[LV2]{LV2}  Lapointe, L. and Vinet, L. : Exact operator solution of
the Calogero-Sutherland model; Comm. Math. Physics 178(1996), 425-452.

\bibitem[Sa]{Sa}  Sahi, S. : A new scalar product for nonsymmetric Jack
polynomials; IMRN 20(1996), 997-1004.

\bibitem[St]{St}  Stanley, R. P. : Some combinatorial properties of Jack
symmetric functions; Adv. Math. 77(1989), 76-115.

\bibitem[Y]{Y}  Yamamoto, T.: Multicomponent Calogero model of $B_N$-type
confined in harmonic potential; Phys. Lett. A208(1995), 293-302.

\bibitem[YT]{YT}  Yamamoto, T. and Tsuchiya, O. : Integrable $1/r^2$ spin
chain with reflecting end; J. Phys. A29(1996), 3977-3984.

\bibitem[Yn]{Yn}  Yan, Z. : A class of hypergeometric functions in several
variables; Canad. J. Math. 44(1992), 1317-1338.
\end{thebibliography}
\end{document}